\documentclass[12pt]{amsart}

\usepackage[scanall]{psfrag}
\usepackage{graphicx}

\newcommand{\bpf}{\begin{proof}}
\newcommand{\epf}{\end{proof}}

\newtheorem{thm}{Theorem}[section]
\newtheorem{cor}[thm]{Corollary}
\newtheorem{lem}[thm]{Lemma}
\newtheorem{prop}[thm]{Proposition}

\newtheorem{case}{Case}

\newtheorem{claim*}{Claim}

\theoremstyle{definition}

\newcommand{\set}[1]{\left\{#1\right\}}

\newcommand{\ben}{\begin{enumerate}}
\newcommand{\een}{\end{enumerate}}
\newcommand{\ita}{\item[(a)]}
\newcommand{\itb}{\item[(b)]}
\newcommand{\itc}{\item[(c)]}
\newcommand{\itd}{\item[(d)]}

\newcommand{\bit}{\begin{itemize}}
\newcommand{\eit}{\end{itemize}}

\newcommand{\wh}{\widehat}
\newcommand{\mc}{\mathcal}

\newcommand{\ove}{\overline{e}}
\newcommand{\bg}{\overline{G}}
\newcommand{\bgt}{\overline{G}_T}
\newcommand{\bgs}{\overline{G}_S}

\newcommand{\cl}{\text{cl}\,}

\newcommand{\Figw}[4]{
\begin{center}
\includegraphics[width=#1]{#2}
\caption{ #3 \label{#4} }
\end{center}}

\begin{document}

\title{Correction to ``Toroidal and Klein bottle boundary slopes''}%

\author[L. G. Valdez-S\'anchez]{Luis G. Valdez-S\'anchez}
\address{Department of Mathematical Sciences,
University of Texas at El Paso\\
El Paso, TX 79968, USA}
\email{valdez@math.utep.edu}%

%

\begin{abstract}
Let $M$ be a compact, connected, orientable, irreducible 3-manifold and
$T_0$ an incompressible torus boundary component of $M$ such that the
pair $(M,T_0)$ is not cabled. In the paper ``Toroidal and Klein bottle
boundary slopes'' \cite{valdez11} by the author it was established that
for any $\mc{K}$-incompressible tori $F_1,F_2$ in $(M,T_0)$ which
intersect in graphs $G_{F_i}=F_i\cap F_j\subset F_i$, $\{i,j\}=\{1,2\}$,
the maximal number of mutually parallel, consecutive, negative edges that
may appear in $G_{F_i}$ is $n_j+1$, where $n_j=|\partial F_j|$. In this
paper we show that the correct such bound is $n_j+2$, give a partial
classification of the pairs $(M,T_0)$ where the bound $n_j+2$ is reached,
and show that if $\Delta(\partial F_1,\partial F_2)\geq 6$ then the bound
$n_j+2$ cannot be reached; this latter fact allows for the short proof of
the classification of the pairs $(M,T_0)$ with $M$ a hyperbolic
3-manifold and $\Delta(\partial F_1,\partial F_2)\geq 6$ to work without
change as outlined in \cite{valdez11}.
\end{abstract}

\maketitle

\section{Introduction}\label{intro}

Let $M$ be a compact, connected, orientable, irreducible 3-manifold and
$T_0$ an incompressible torus boundary component of $M$ such that the
pair $(M,T_0)$ is not cabled and {\it irreducible} (that is, $M$ is
irreducible and $T_0$ is incompressible in $M$). A punctured torus
$(F,\partial F)\subset (M,T_0)$ is said to be {\it generated by a (an
essential) Klein bottle} if there is a (an essential, resp.) punctured
Klein bottle $(P,\partial P)\subset (M,T_0)$ such that $F$ is isotopic in
$M$ to the frontier of a regular neighborhood of $P$ in $M$. We also say
that $F$ is {\it $\mc{K}$-incompressible} if $F$ is either incompressible
or generated by an essential Klein bottle.

The main purpose of the paper \cite{valdez11} was to establish an upper
bound for the maximal number of mutually parallel, consecutive, negative
edges that may appear in either graph of intersection $G_{F_1},G_{F_2}$
between $\mc{K}$-incompressible punctured tori $F_1,F_2$ in $(M,T_0)$. In
\cite[Proposition 3.4]{valdez11} it is proved that for $\{i,j\}=\{1,2\}$
and $n_j=|\partial F_j|$, if $G_{F_i}$ contains such a collection of
$n_j+2$ negative edges then $M$ is homeomorphic to the trefoil knot
exterior or to one of the manifolds $P\times S^1/[m]$, $m\geq 1$
constructed in \cite[\S 3.4]{valdez11}, none of which is a hyperbolic
manifold; consequently, if $M$ is not one of the manifolds listed in
\cite[Proposition 3.4]{valdez11}, the upper bound for such a
collection of negative edges was found to be $n_j+1$.

In this paper we show that the list of options for the homeomorphism
class of $M$ given in \cite[Proposition 3.4]{valdez11} is incomplete, so
that if $M$ is not one of the manifolds listed in \cite[Proposition
3.4]{valdez11} then the correct bound for such families of negative edges
in the graph $G_{F_i}$ is $n_j+2$, and that if the upper bound $n_j+2$ is
reached then $(M,T_0)$ belongs to a certain family of examples each of
which contains a separating essential twice punctured torus with boundary
slope at distance 3 from that of $F_j$.

We will use the same notation set up in \cite{valdez11} except that {\it
polarized} will be replaced by {\it positive} (see Section~\ref{pre}); in
particular, the tori $F_1,F_2$ will now be denoted by $S,T$, with
$s=|\partial S|$ and $t=|\partial T|$. A graph $G$ in a punctured surface
$F$ is a 1-submanifold properly embedded in $F$ with vertices the
components of $\partial F$. The graph $G$ is {\it essential} if no edge
is parallel into $\partial F$ and each circle component is essential in
$F$. The reduced graph $\bg\subset F$ of $G$ is the graph obtained from
$G$ by amalgamating each maximal collection of mutually parallel edges
$e_1,\dots,e_k$ of $G$ into a single arc $\ove\subset F$; we then say
that the {\it size} of $\ove$ is $|\ove|=k$. The symbol
$(+g,b\,;\alpha_1/\beta_1,\dots,\alpha_k/\beta_k)$ will be used to denote
a Seifert fibered manifold over an orientable surface of genus $g\geq 0$
with $b\geq 0$ boundary components and $k\geq 0$ singular fibers of
orders $\beta_1,\dots,\beta_k$.

The main technical result of this paper is the following.

\begin{prop}\label{main}
For each $t\geq 4$ there is an irreducible pair $(M_t,T_0)$ with
$\partial M_t=T_0$ which contains properly embedded essential punctured
tori $(S,\partial S)$ and $(T,\partial T)\subset (M,T_0)$ satisfying the
following properties:
\ben
\item
$S$ is a separating twice punctured torus and $T$ is a positive torus
with $|\partial T|=t$,

\item
$\Delta(\partial S,\partial T)=3$,

\item
$S,T$ intersect transversely and minimally in the essential graphs
$G_T=S\cap T\subset T$ shown in Figs.~\ref{n21-4} and
\ref{n21-5}, where the reduced graph $\bgs$ consists of 3
negative edges of sizes $t+2,t,t-2$,

\item $M_t(\partial T)$ is a torus bundle over the circle with fiber
$\wh{T}$, and $M_t$ is not homeomorphic to any of the manifolds $P\times
S^1/[m]$ constructed in \cite[Proposition 3.4]{valdez11}.
\een
\end{prop}

It is proved in \cite{valdez13} (preprint, in progress) that each
manifold $M_t$ in Proposition~\ref{main} is hyperbolic with $M_t(\partial
S)=(+0,1:\alpha_1/2,\alpha_2/(t+2))\cup_{\wh{S}}
(+0,1:\alpha_1/2,\alpha_2/(t-2))$, so $M_{t_1}\not\approx M_{t_2}$ for
$t_1\neq t_2$ and $S$ is generated by a punctured Klein bottle iff $t=4$.

The corrected version of \cite[Proposition 3.4]{valdez11} can now be
stated as follows.

\begin{prop}(Correction to \cite[Proposition 3.4]{valdez11})
\label{prop2}
Let $(M,T_0)$ be an irreducible pair which is not cabled, $(T,\partial
T)\subset (M,T_0)$ a $\mc{K}$-incompressible torus with $t=|\partial
T|\geq 1$, and $R\subset M$ a surface which intersects $T$ in essential
graphs $G_R,G_T$, such that $G_R$ has at least $t+2$ mutually parallel,
consecutive negative edges. Then $T$ is a positive torus and one of the
following holds:

\ben
\item
the conclusion of \cite[Proposition 3.4]{valdez11} holds, so $M$ is
homeomorphic to the trefoil knot exterior $(t=1)$ or to one of the
manifolds $P\times S^1/[m]$ $(t\geq 2)$,

\item
$t\geq 4$ and $(M,T_0)$ is homeomorphic to one of the pairs $(M_t,T_0)$
of Proposition~\ref{main}, and $|\ove|\leq t+2$ holds for any edge of the
reduced graph $\bg_R$,

\item
$t=2$ with $M=(+0,1;-1/4,-1/4)$ and $M(\partial T)=(+0,0;1/2,-1/4,-1/4)$,
or $t=3$ with $M=(+0,1;-1/3,-1/6)$ and $M(\partial
T)=(+0,0;1/2,-1/3,-1/6)$, where in each case the essential annulus
$A\subset (M,T_0)$ satisfies $\Delta(\partial T,\partial A)=2$.
\een
\end{prop}

The smaller bound of $n_j+1$ allows for the short proof  of the
classification of hyperbolic manifolds $(M,T_0)$ with toroidal or
Kleininan Dehn fillings at distance $6\leq\Delta\leq 8$ given in \cite[\S
4]{valdez11}. The following result states that in the range
$6\leq\Delta\leq 8$ the bound $n_j+2$ is never reached, which implies
that the proofs in \cite[\S 4]{valdez11} work as written.

\begin{lem}\label{corr2}
Suppose $(M,T_0)$ is a hyperbolic manifold and $S,T\subset (M,T_0)$ are
essential tori such that $t=|\partial T|\geq 3$ and $\Delta(\partial
S,\partial T)\geq 6$; then $|\ove|\leq t+1$ holds in $\bgs$.
\end{lem}

Our last lemma summarizes the bounds on the sizes of negative edges of
reduced graphs like $\bg_R$ obtained from the above results.

\begin{lem}\label{tp2a}
Let $(M,T_0)$ be an irreducible pair which is not cabled, $(T,\partial
T)\subset (M,T_0)$ a $\mc{K}$-incompressible torus with $t=|\partial
T|\geq 1$, and $Q\subset M$ a surface which intersects $T$ in essential
graphs $G_Q=Q\cap T\subset Q,G_T=Q\cap T\subset T$. Then one of the
following holds:
\ben
\item
$t\leq 3$ and $M$ is one of the Seifert manifolds $(+0,1;1/2,1/3)$
$(t=1)$, $(+0,1;-1/4,-1/4)$ $(t=2)$, or $(+0,1;-1/3,-1/6)$ $(t=3)$,

\item
$t\geq 2$ and $M$ is one of the manifolds $P\times S^1/[m]$ constructed
in \cite[Proposition 3.4]{valdez11},

\item
$t\geq 4$ and $M$ is homeomorphic to one of the manifolds $M_t$ of
Proposition~\ref{main}, in which case the bound $|\ove|\leq t+2$ holds
for all negative edges of $\bg_Q$,

\item
the bound $|\ove|\leq t+1$ holds for all negative edges of $\bg_Q$.
\een
\end{lem}

In Section~\ref{slide} we review the notation and constructions given in
\cite[\S 2,\S 3]{valdez11}, with which we assume the reader is familiar.
In Section~\ref{t3} we construct the manifolds $M_t$ for $t\geq 4$ of
Proposition~\ref{main} and establish the results needed to prove
Proposition~\ref{prop2} and Lemmas~\ref{corr2} and \ref{tp2a}.

\section{Slidable and non-slidable bigons}\label{slide}

\subsection{Generalities}\label{pre}

For any slope $r$ in $T_0\subset\partial M$, $M(r)$ denotes the Dehn
filling $M\cup_{T_0}S^1\times D^2$, where $r$ bounds a disk in the solid
torus $S^1\times D^2$. We denote the core of $S^1\times D^2$ by
$K_r\subset M(r)$.

Let $F$ be a surface properly embedded in $(M,T_0)$ with $|\partial
F|\geq 1$ and boundary slope $r$. Then the surface $\wh{F}\subset M(r)$
obtained by capping off the components of $\partial F$ with a disjoint
collection of disks in $S^1\times D^2$ is a closed surface, which we
always assume to intersect $K_r$ transversely and minimally in
$M(\partial F)$.

If $F$ is orientable, we say that $F$ is {\it neutral} if $\wh{F}\cdot
K_r=0$ in $M(\partial F)$, where $\wh{F}\cdot K$ denotes homological
intersection number, and that $F$ is {\it positive} if $|\wh{F}\cdot
K_r|=|\wh{F}\cap K_r|$; the latter is equivalent to the term {\it
polarized} used in
\cite{valdez11}.

For surfaces $F_1,F_2\subset(M,T_0)$ with transverse intersection, for
$i=1,2$, $G_{F_i}=F_1\cap F_2\subset F_i$ will denote the graph of
intersection of $F_1$ and $F_2$ in $F_i$ (with vertices the boundary
components of $F_i$).

Following \cite{tera10}, for $i=1,2$ we orient the components of
$\partial F_i$ and coherently on $T_0$ and say that an edge $e$ of
$G_{F_i}$ is {\it positive} or {\it negative} depending on whether the
orientations of the components of $\partial F_i$ (possibly the same)
around a small rectangular regular neighborhood of $e$ in $F_i$ appear as
in Fig~\ref{n20}.
\begin{figure}
\psfrag{e}{$e$}
\psfrag{positive edge}{positive edge}
\psfrag{negative edge}{negative edge}
\Figw{2.5in}{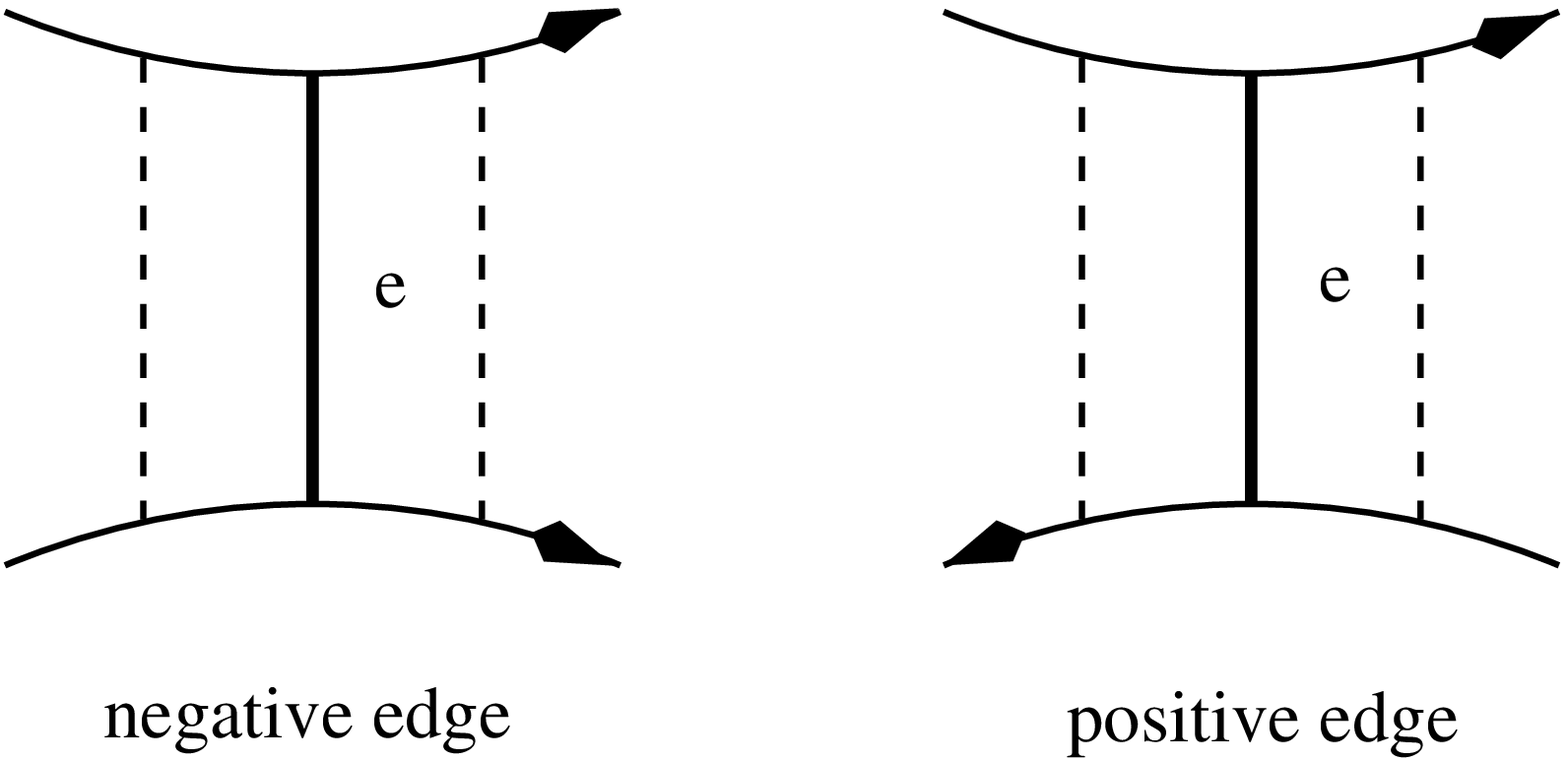}{}{n20}
\end{figure}

The following lemma summarizes some of the general properties of graphs
of intersection of surfaces in $(M,T_0)$ that will be relevant in the
sequel.

\begin{lem}\label{basic}
Let $F_1,F_2$ be properly embedded surfaces in $(M,T_0)$ with essential
graphs of intersection $G_{F_{1}}=F_1\cap F_2\subset F_1$ and
$G_{F_{2}}=F_1\cap F_2\subset F_2$.
\ben
\ita
Parity Rule: for $\{i,j\}=\{1,2\}$, an edge of $F_1\cap F_2$ is positive
in $G_{F_{i}}$ iff it is negative in $G_{F_{j}}$ (cf.\ \cite{tera10});
\een

moreover, if $(M,T_0)$ is not cabled,

\ben
\itb
no two edges of $F_1\cap F_2$ are parallel in both $G_{F_1}$ and
$G_{F_2}$ (\cite[Lemma 2.1]{gordon5}),

\itc
if $F_i$ is a torus and $G_{F_{j}}$ has a family of $n_{i}+1$ mutually
parallel, consecutive, negative edges then $F_i$ is a positive torus
(cf.\ \cite[Lemma 3.2]{valdez11}).
\hfill\qed
\een
\end{lem}

\subsection{Review of constructions in
\cite[\S 2,3]{valdez11}.}
Suppose $(M,T_0),T,R$ satisfy the hypothesis in Proposition~\ref{prop2},
so $(M,T_0)$ is an irreducible pair with $T_0$ a torus boundary component
of $M$ which is not cabled, $T$ a $\mc{K}$-incompressible torus in
$(M,T_0)$, and $R$ any surface properly embedded in $(M,T_0)$ which
intersects $T$ transversely in essential graphs $G_T=R\cap T\subset T$
and $G_R=R\cap T\subset R$.

We assume there is a collection $E=\{e_1,e_2,\dots,e_{t+1},e_{t+2}\}$ of
mutually parallel, consecutive, negative edges in $G_R$. By
Lemma~\ref{basic}(c), the torus $T$ is actually positive and hence
incompressible in $M$.

The torus $T$ has orientation vector $\vec{N}$ shown in Fig.~\ref{n04},
and each vertex of $T$ is given the orientation induced by $\vec{N}$ (see
Fig.~\ref{n04}). Recall that the vertices of $G_T,G_R$ are the components
of $\partial T,\partial R$, respectively; we denote the vertices of $T$
by $v_i$'s and those of $R$ by $w_k$'s. Any two edges $e,e'$ of $G_T$
that are incident to the oriented vertex $v$ split $v$ into two open
subintervals $(e,e')\subset v$ and $(e',e)\subset v$, where $(e,e)$ is
the open subinterval whose {\it left} and {\it right} endpoints, as
defined by the orientation of $v$, come from $e$ and $e'$, respectively.

The collection of edges $E$ induces a permutation $\sigma$ of the form
$x\mapsto x+\alpha$ with $1\leq
\alpha\leq t$, where $\gcd(t,\alpha)=1$ by
\cite[Lemma 3.2]{valdez11}; the definition of $\sigma$
requires that a common orientation be given to the edges of $E$, and
reversing the orientation of such edges replaces $\sigma$ with its
inverse, hence $\alpha$ with $t-\alpha$.

Cutting $M$ along $T$ produces an irreducible manifold $M_T=\cl
(M\setminus N(T))$ with copies $T^1,T^2\subset\partial M_T$ of $T$ on its
boundary and strings $I'_{1,2},I'_{2,3},\dots,I'_{t,1}\subset\partial
M_T$ such that $M_T/\psi=M$ for some orientation preserving homeomorphism
$\psi:T^1\to T^2$, where $T^1,T^2$ are oriented by normal vectors
$\vec{N}^1,\vec{N}^2$ as shown in Fig.~\ref{n05-2}.

We assume that the edges of $E$ and $T^1\cap R,T^2\cap R$ are arranged in
$G_R$ as shown in Fig.~\ref{n11}. Each edge $e_i$ of $E$ and vertex $v_i$
of $T$ gives rise to two copies of itself $e^1_i, v^1_i\subset T^1$ and
$e^2_i,v^1_i\subset T^2$ such that $\psi(e^1_i)=e^2_i$ and
$\psi(v^1_i)=v^2_i$.

The collection $E$ also gives rise to two essential cycles in $T$,
$\gamma_1=e_1\cup e_2\cup\cdots\cup e_t$ and $\gamma_2=e_2\cup
e_3\cup\cdots\cup e_t\cup e_{t+1}$, such that
$\Delta(\gamma_1,\gamma_2)=1$ holds in $\wh{T}$. The edges $e^1_1\cup
e^1_2\cup\cdots\cup e^1_t$ form the essential cycle $\gamma^1_1$ in
$T^1$, while $e^2_2\cup e^2_3\cup\cdots\cup e^2_{t+1}$ form the essential
cycle $\gamma^2_2$ in $T^2$, such that the bigon faces $F'_1,F'_2,\dots,
F'_t$ bounded by $e_1,e_2,\dots, e_{t+1}$ in $R\cap M_T$ form an
essential annulus $\mc{A}$ in $M_T(\partial T)$ as shown in
Fig.~\ref{n05-2}. For simplicity, we refer to the union of the bigons
$F'_1,F'_2,\dots,F'_t$ in $M_T$ also as the annulus $\mc{A}$.

The next result describes the embedding of $\mc{A}\cup F'_{t+1}$ in $M_T$
and the structures of $M_T$ and $M_T(\partial T)$; its proof follows
immediately from the arguments of \cite[Lemma 3.3]{valdez11}.

\begin{lem}{(\cite[Lemma 3.3]{valdez11})}\label{mt}
Up to homeomorphism, the bigons $F'_1,\cdots,F'_{t+1}$ lie in $M_T$ as
shown in Fig.~\ref{n05-2}. In particular, $M_T\approx T\times I$ is a
genus $t+1$ handlebody with $F'_1,\cdots,F'_{t+1}$ a complete disk
system, $\partial M=T_0$, and $M(\partial T)$ is a torus bundle over the
circle with fiber $\wh{T}$.\hfill\qed
\end{lem}

\begin{figure}
\psfrag{ui}{$w_i$}
\psfrag{ui'}{$w_{j}$}
\psfrag{e1}{$e_1$}
\psfrag{e2}{$e_2$}
\psfrag{e3}{$e_3$}
\psfrag{et}{$e_t$}
\psfrag{et1}{$e_{t+1}$}
\psfrag{t1}{$T^1$}
\psfrag{t2}{$T^2$}
\psfrag{f1}{$F'_1$}
\psfrag{f2}{$F'_2$}
\psfrag{ft}{$F'_t$}\psfrag{ft1}{$F'_{t+1}$}
\psfrag{e21}{$e^2_1$}
\psfrag{e11}{$e^1_1$}
\psfrag{e2t1}{$e^2_{t+1}$}
\psfrag{e1t1}{$e^1_{t+1}$}
\psfrag{e2t2}{$e^2_{t+2}$}
\psfrag{e1t2}{$e^1_{t+2}$}
\psfrag{I'12}{$I'_{1,2}$}
\psfrag{I't1}{$I'_{t,1}$}
\Figw{5.5in}{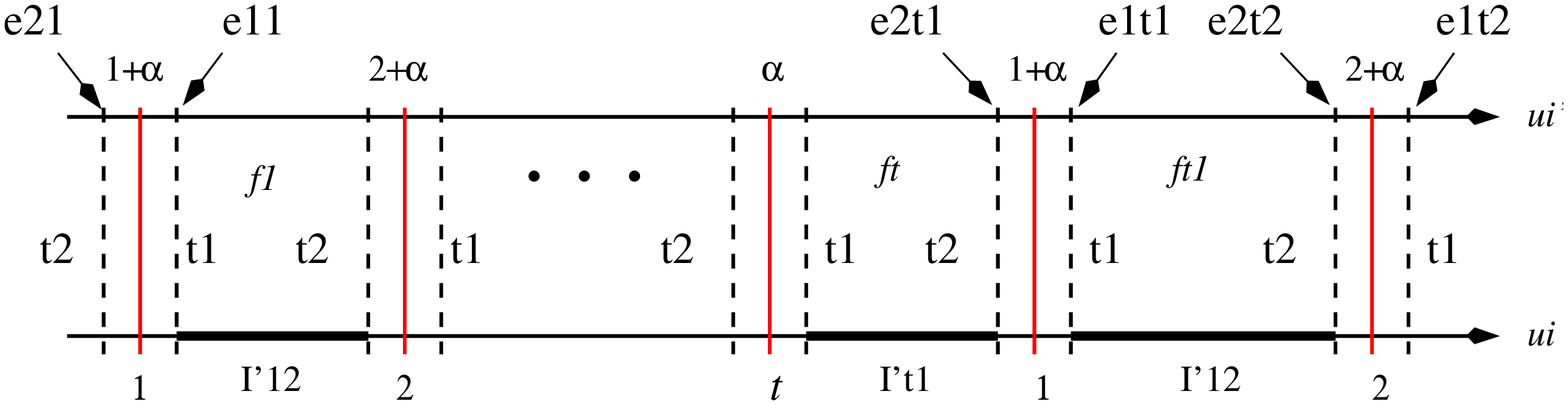}{}{n11}
\end{figure}

\begin{figure}
\psfrag{g}{$\gamma_1$}
\psfrag{e1}{$e_1$}
\psfrag{et1}{$e_{t+1}$}
\psfrag{mt}{$\mu_{t+1-\alpha}$}
\psfrag{m1}{$\mu_{1}$}
\psfrag{m1a}{$\mu_{1+\alpha}$}
\psfrag{N}{$\vec{N}$}\psfrag{T}{$T$}
\Figw{3.5in}{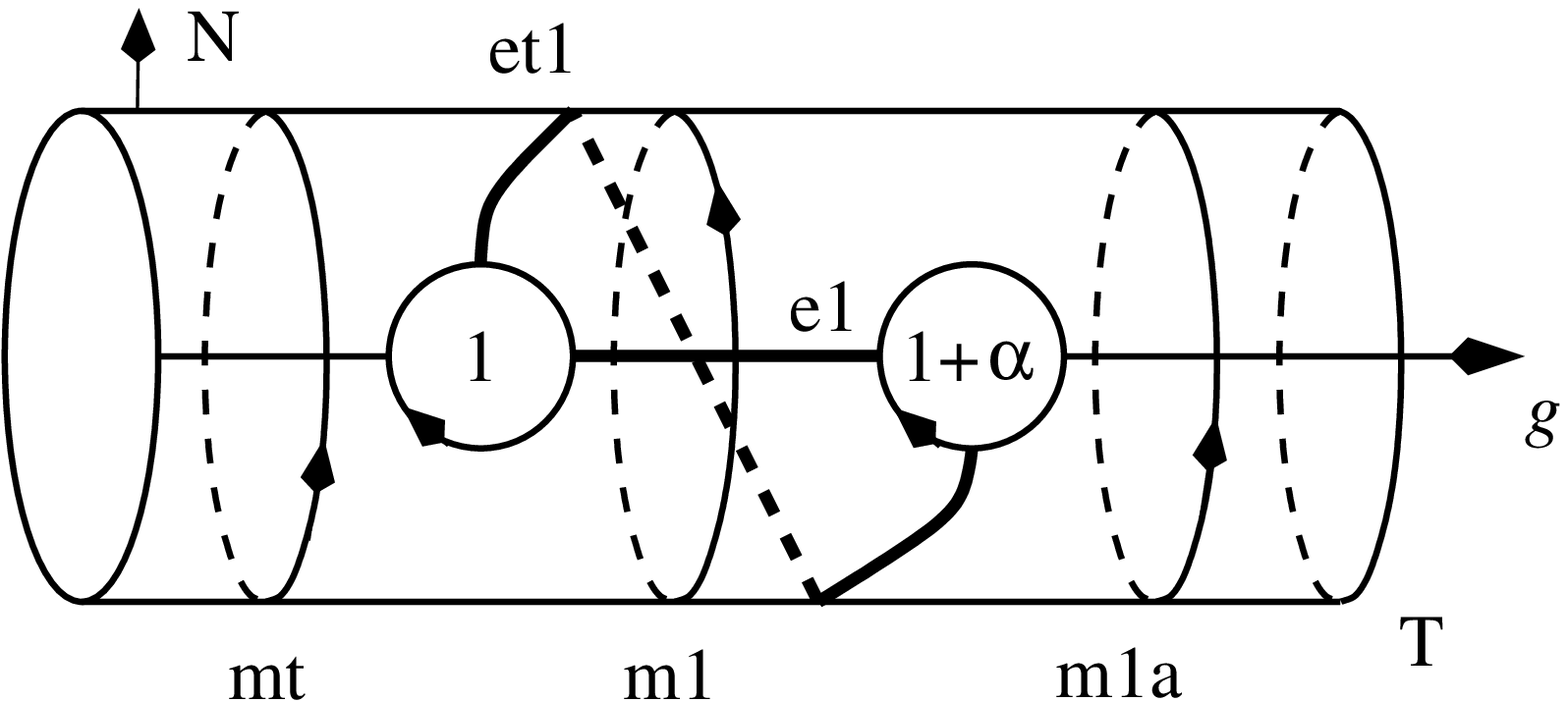}{}{n04}
\end{figure}

\begin{figure}
\psfrag{g1}{$\gamma_1^1$}
\psfrag{g2}{$\gamma_2^2$}
\psfrag{e22}{$e^2_2$}
\psfrag{e11}{$e^1_1$}
\psfrag{e21}{$e^2_1$}
\psfrag{u11}{$\mu^1_1$}
\psfrag{u22}{$\mu^2_2$}
\psfrag{u1k}{$\mu^1_k$}
\psfrag{u2k1}{$\mu^2_{k+1}$}
\psfrag{N1}{$\vec{N}^1$}
\psfrag{N2}{$\vec{N}^2$}
\psfrag{e2t2}{$e^2_{t+2}$}
\psfrag{e1t1}{$e^1_{t+1}$}
\psfrag{e1t}{$e^1_{k}$}
\psfrag{T1}{$T^1$}
\psfrag{T2}{$T^2$}
\Figw{5in}{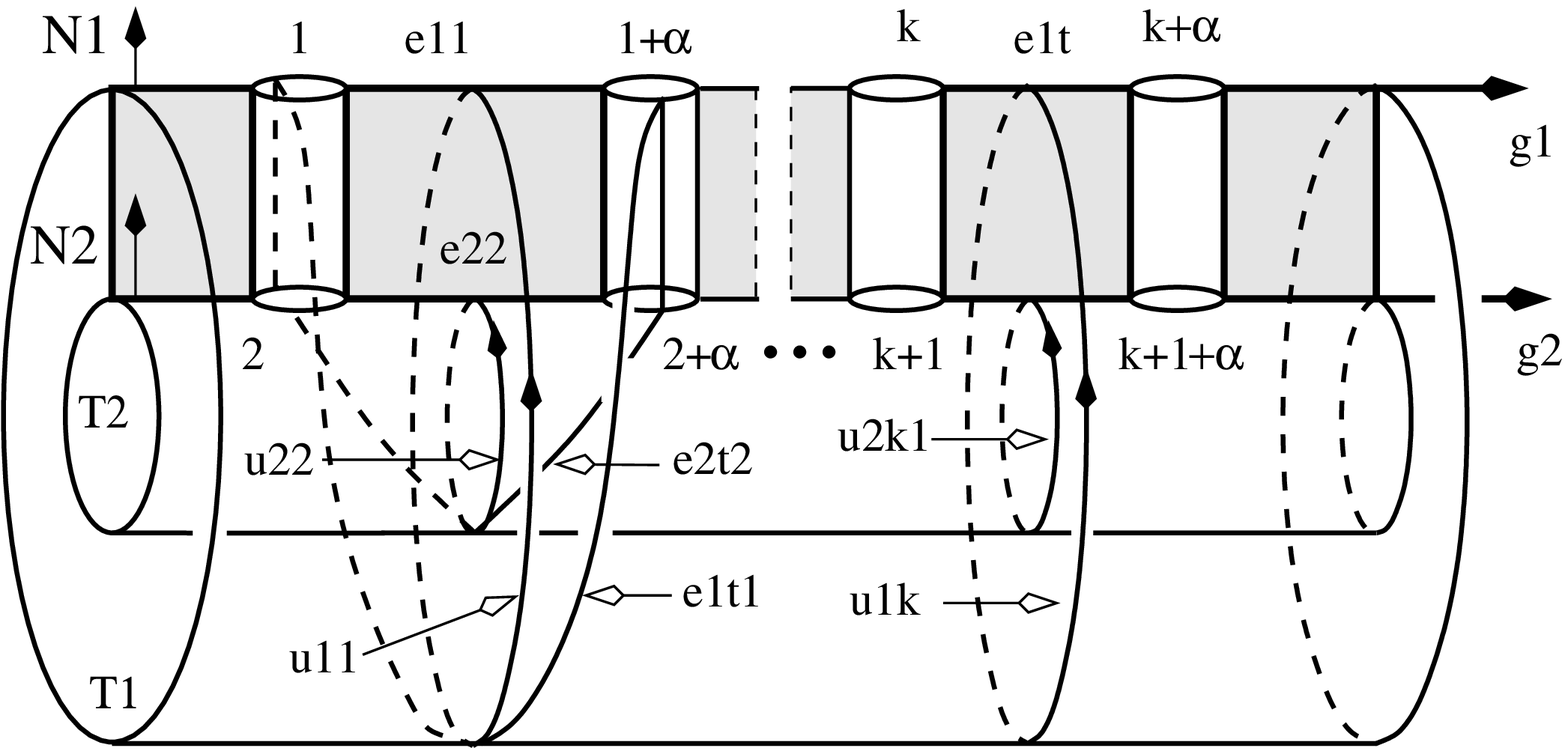}{}{n05-2}
\end{figure}

Finally, we construct auxiliary circles $\mu_i$ in $T$ having the same
slope in $\wh{T}$ as the cycle $e_1\cup e_{t+1}$, oriented and labeled as
shown in Fig.~\ref{n04}. The counterparts $\mu^1_i\subset T^1$ of the
$\mu_i$'s are shown in Fig.~\ref{n05-2}, while the circles
$\mu^2_i\subset T^2$ are represented abstractly in Fig.~\ref{n05-2} since
the location of the edge $e^2_1\subset T^2$ is not given yet.

\subsection{Review of the argument of \cite[Proposition
3.4]{valdez11}} At this point the argument used in the proof of
\cite[Lemma 3.6]{valdez11} states that, for $t\geq 2$,

\begin{quote}
{\it \dots the faces $F'_1$ and $F'_{t+1}$ can be isotoped in $M_T$ to
construct an annulus $A_1\subset M_T$ with boundary the circles
$\mu^1_1\cup\mu^2_2$, which under their given orientations remain
coherently oriented relative to $A_1$. Via the product structure
$M_T=T\times I$, it is not hard to see that each pair of circles
$\mu^1_k,\mu^2_{k+1}$ cobounds such an annulus $A_k\subset M_T$ for
$1\leq k\leq t$, with the oriented circles $\mu^1_k,\mu^2_{k+1}$
coherently oriented relative to $A_k$; these annuli $A_k$ can be taken to
be mutually disjoint and $I$-fibered in $M_T=T\times I$. Since
$\psi(\mu^1_k)=\mu^2_k$ (preserving orientations), the union $A_1\cup
A_2\dots\cup A_t$ yields a closed nonseparating torus $T''$ in $M$, on
which the circles $\mu_1,\mu_2,\dots,\mu_t$ appear consecutively in this
order and coherently oriented.}
\end{quote}

The problem with the above argument is that it is assumed from the
beginning that the boundary of the annulus $A_1$ must necessarily be
$\mu^1_1\cup\mu^2_2$, which, as we shall see next, is not the case and
leads to the present correction to
\cite[Proposition 3.4]{valdez11}.

The isotopy of $F'_1$ and $F'_{t+1}$ in $M_T$ mentioned in the above
quote can be thought of as the result of a {\it sliding process}, where
the corners of the face $F'_{t+1}$ are slid onto the face $F'_1$ so as to
coincide with each other, at which point the isotoped face $F'_{t+1}$
becomes the annulus $A_1$ properly embedded in $M_T$ with boundary the
circles $\mu^1_1\cup\mu^2_2$.

We will say that the bigon $F'_{t+1}$ is {\it slidable (relative to the
annulus $\mc{A}$)} whenever the annulus $A_1$ produced by the above
isotopy of $F'_{t+1}$ satisfies $\partial A_1=\mu^1_1\cup\mu^2_2$, and
otherwise that $F'_{t+1}$ is {\it non-slidable}. Equivalently, $F'_{t+1}$
is slidable iff the cycles $\psi(e^1_1\cup e^1_{t+1})=e^2_1\cup
e^2_{t+1}\subset T^2$ and $e^2_2\cup e^2_{t+2}\subset T^2$ have the same
slope in $\wh{T}^2$, that is iff the cycles $e_1\cup e_{t+1}$ and
$e_2\cup e_{t+2}$ have the same slope in $\wh{T}$.

As we shall see in the next section, there are two combinatorially
different embeddings of the edge $e^2_1$ in $T^2$ which correspond to the
annulus $A_1$ being slidable or not; the generic embeddings $e^2_1\subset
T^2$ are shown in Fig.~\ref{n21-2a}, the slidable case which produces the
circles $\mu^2_i\subset T^2$ shown in Fig.~\ref{n05-2}), and
Fig.~\ref{n21-2b}, the non-slidable case.

Using this notation we summarize \cite[Proposition 3.4]{valdez11} as
follows:

\begin{lem}{\cite[Proposition
3.4]{valdez11}}\label{lem01}
\ben
\item
If $t=1$ then $M$ is the exterior of the trefoil knot $(+0,1;1/2,1/3)$.

\item
If $t\geq 2$ and $F'_{t+1}$ is slidable then $M$ is homeomorphic to one
of the manifolds $P\times S^1/[m]$ constructed in
\cite[\S 3.4]{valdez11}, in which case $M$ is not Seifert fibered
and contains a closed nonseparating torus.\hfill\qed
\een
\end{lem}

We shall see below that in most cases, which include those with $t\geq 4$
and $\alpha\not\equiv \pm 1\mod t$, the bigon $F'_{t+1}$ is slidable, and
that the exceptions with $t\geq 4$ form the family $M_t$ of
Proposition~\ref{main}.

\section{Main results}\label{t3}

In this section we assume that $(M,T_0),T,R$ satisfy the hypothesis of
Proposition~\ref{prop2}; also, the mutually parallel edges
$E=\{e_1,e_2,\dots,e_{t+1},e_{t+2}\}$ in $G_R$ are labeled as in
Fig.~\ref{n11}, induce the permutation $\sigma(x)\equiv x+\alpha\mod t$
with $\gcd(t,\alpha)=1$, and cobound bigon faces $F'_1,\dots,F'_{t+1}$
embedded in $M_T$ as shown in Fig.~\ref{n05-2}. The case $t=1$ is
considered in Lemma~\ref{lem01}(1).

\subsection{The cases $t=2,3$}

\begin{lem}\label{t2t3}
If $t=2,3$ then $(M,T_0)$ satisfies the conclusion of
\cite[Proposition 3.4]{valdez11} or of
Proposition~\ref{prop2}(3).
\end{lem}

\bpf
For $t=2$ we must have $\alpha=1$; for $t=3$ we may also assume that
$\alpha=1$ after reversing the orientation of the edges of $E$ if
necessary.

We begin with a detailed analysis of the case $t=2$. Fig.~\ref{n21c}(a)
shows the edges $e_i$ and bigons $F'_i$ of $E$ in $G_R$. By
Lemma~\ref{mt}, up to homeomorphism, the annulus $\mc{A}$ cobounded by
the cycles $\gamma_1^1=e^1_1\cup e^1_2\subset T^1$ and
$\gamma_2^2=e^2_2\cup e^2_3\subset T^2$ and the bigon $F'_3$ cobounded by
$e_3,e_4$ may be assumed to lie in $M_T$ as shown in Fig.~\ref{n21c}(b)
or (c); for simplicity, the upper labels in the edges will not be shown
in the figures representing $M_T$.

Notice that the embeddings of the edges $e_2$ and $e_3$ are determined in
both $T^1$ and $T^2$ at this point, but that the embeddings of $e_1$ and
$e_4$ are so far determined only in $T^1$ and $T^2$, respectively.

\begin{figure}
\psfrag{g1}{$\gamma_1^1$}
\psfrag{g2}{$\gamma_2^2$}
\psfrag{e11}{$e^1_1$}
\psfrag{e21}{$e^2_1$}
\psfrag{e12}{$e^1_2$}
\psfrag{e22}{$e^2_2$}
\psfrag{e13}{$e^1_3$}
\psfrag{e23}{$e^2_3$}
\psfrag{e14}{$e^1_4$}
\psfrag{e24}{$e^2_4$}
\psfrag{vj}{$w_i$}\psfrag{vj'}{$w_j$}
\psfrag{F'1}{$F'_1$}\psfrag{F'2}{$F'_2$}\psfrag{F'3}{$F'_3$}
\psfrag{(a)}{$(a)$}\psfrag{(b)}{$(b)$}\psfrag{(c)}{$(c)$}
\psfrag{T1}{$T^1$}
\psfrag{T2}{$T^2$}
\psfrag{e1}{$e_1$}\psfrag{e2}{$e_2$}
\psfrag{e3}{$e_3$}\psfrag{e4}{$e_4$}
\Figw{5in}{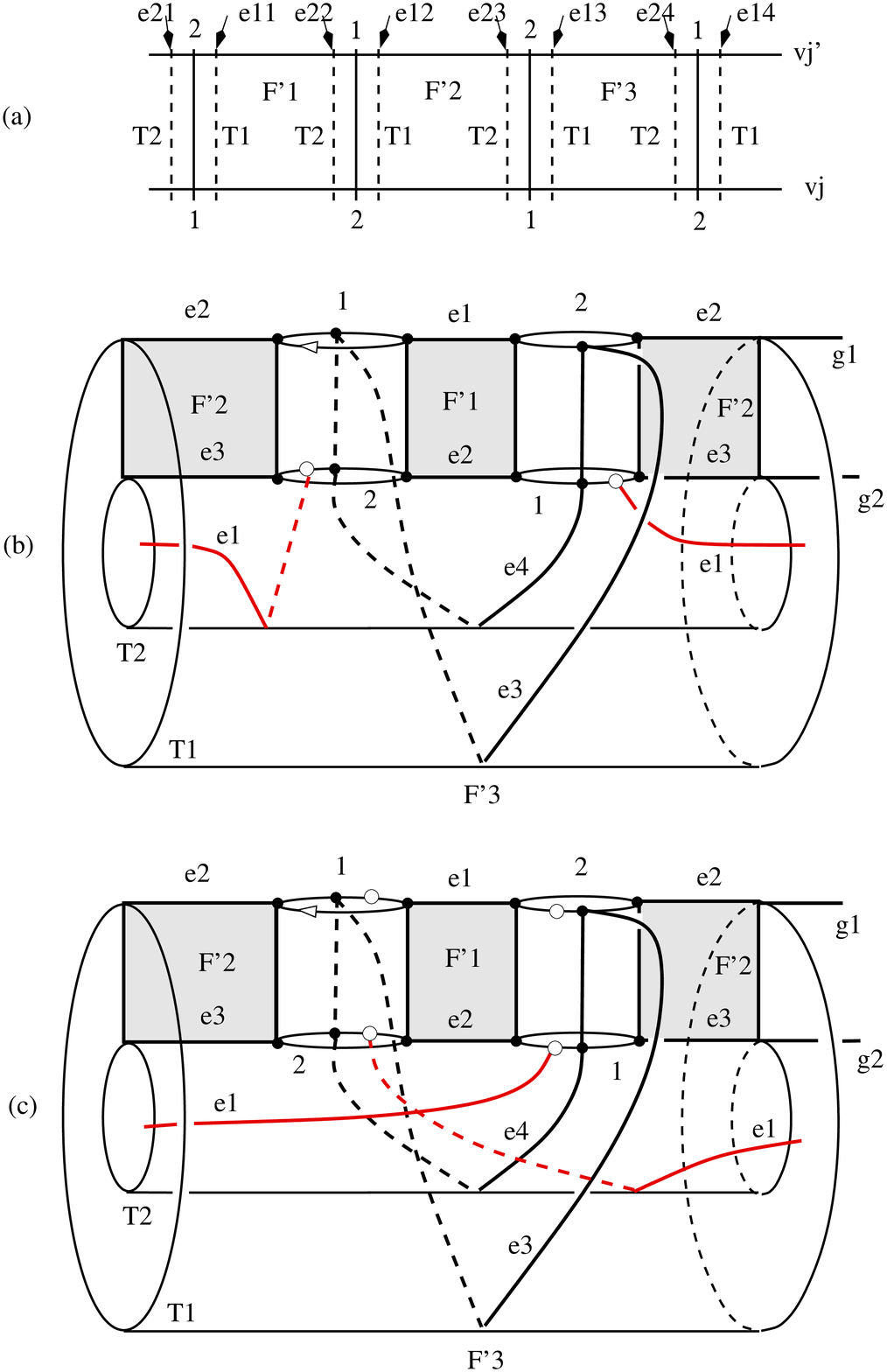}{}{n21c}
\end{figure}

To determine the possible embeddings of $e^2_1$ in $T^2$ we observe that
in either Fig.~\ref{n21c}(b),(c) the endpoint of $e^1_1$ in $v^1_1$ lies
in the interval $(e^1_3,e^1_2)$, and that the same statement holds for
the endpoint of $e^1_1$ in $v^1_2$. Since $\psi$ maps each $v^1_i\subset
T^1$ to $v^2_i\subset T^2$ and each $e^1_i\subset T^1$ to $e^2_i\subset
T^2$, the endpoints of $e^2_1$ in $v^2_1,v^2_2$ must also lie in the
corresponding intervals $(e^2_3,e^2_2)$. As the edge $e^2_4$ is already
embedded in $T^2$ and its endpoint on $v^2_2$ lies in $(e^2_3,e^2_2)$, it
follows that the endpoint of $e^2_1$ in $v^2_2$ must lie in one of the
intervals $(e^2_3,e^2_4)$ or $(e^2_4,e^2_2)$.

The first option is represented in Fig.~\ref{n21c}(b); the placement of
$e^2_1$ in $T^2$ is then uniquely determined using the fact that no two
edges of $E$ are mutually parallel in $T$. Since the cycles $e^2_1\cup
e^2_3$ and $e^2_2\cup e^2_4$ have the same slope in $\wh{T}^2$, the bigon
$F'_{3}$ is slidable and so by Lemma~\ref{lem01} the pair $(M,T_0)$
satisfies the conclusion of \cite[Proposition 3.4]{valdez11}. The
embedding of $e^1_4$ in $T^1$ is now uniquely determined as well: its
endpoints lie in the intervals $(e^1_1,e^1_2)\subset v^1_1$ and
$(e^1_1,e^1_2)\subset v^1_2$.

\begin{figure}
\psfrag{g1}{$\gamma_1^1$}
\psfrag{g2}{$\gamma_2^2$}\psfrag{g}{$\gamma_1$}
\psfrag{e11}{$e^1_1$}
\psfrag{e21}{$e^2_1$}
\psfrag{e12}{$e^1_2$}
\psfrag{e22}{$e^2_2$}
\psfrag{e13}{$e^1_3$}
\psfrag{e23}{$e^2_3$}
\psfrag{e14}{$e^1_4$}
\psfrag{e24}{$e^2_4$}
\psfrag{vj}{$w_i$}\psfrag{vj'}{$w_j$}
\psfrag{F'1}{$F'_1$}\psfrag{F'2}{$F'_2$}\psfrag{F'3}{$F'_3$}
\psfrag{(a)}{$(a)$}\psfrag{(b)}{$(b)$}\psfrag{(c)}{$(c)$}
\psfrag{T1}{$T^1$}
\psfrag{T2}{$T^2$}\psfrag{T}{$T$}
\psfrag{e1}{$e_1$}\psfrag{e2}{$e_2$}
\psfrag{e3}{$e_3$}\psfrag{e4}{$e_4$}
\psfrag{c2}{$c_2$}
\Figw{5in}{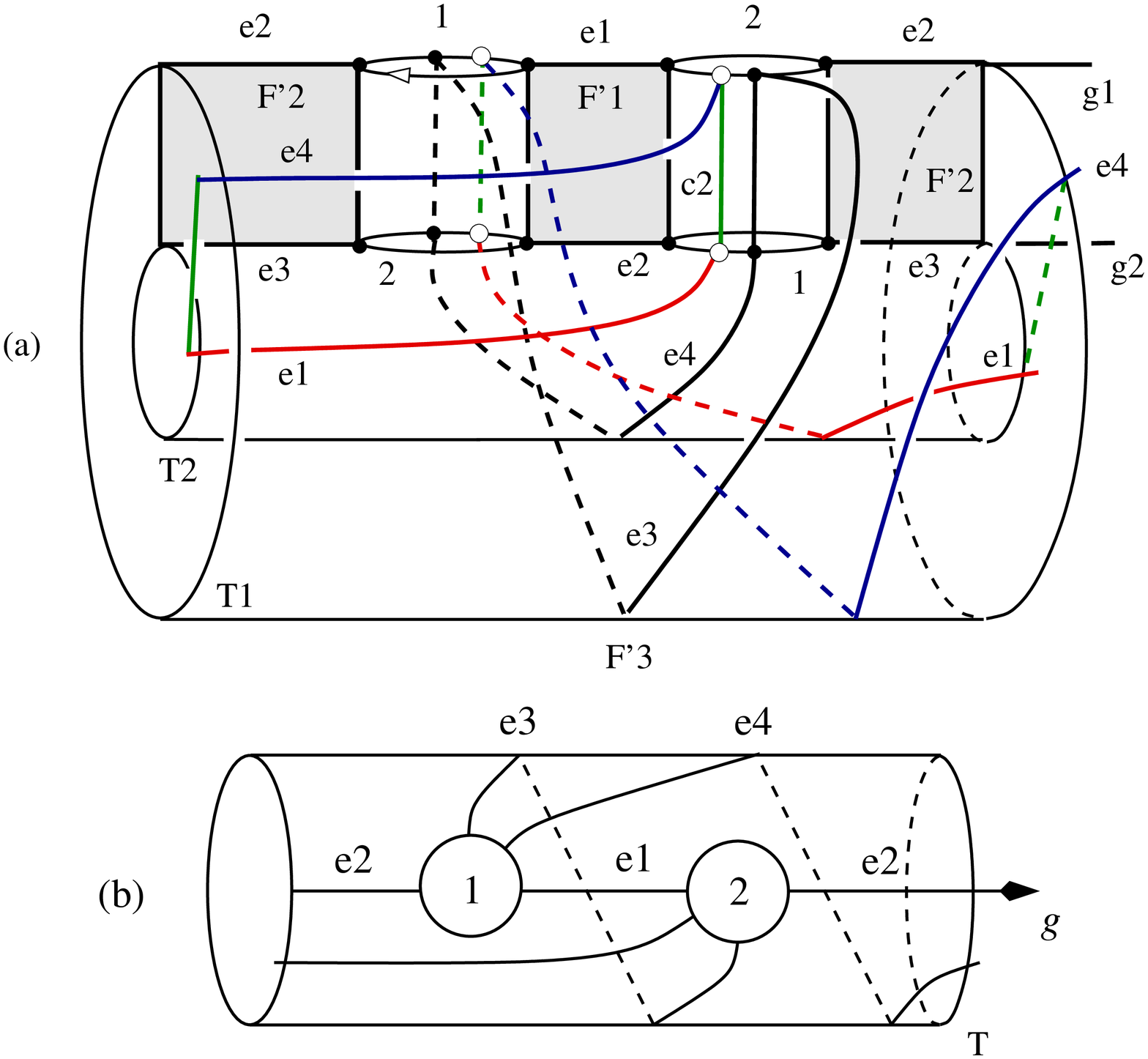}{}{n21c2}
\end{figure}

The second option is represented in Fig.~\ref{n21c}(c), and here we have
that $\Delta(e^2_1\cup e^2_3,e^2_2\cup e^2_4)=2\neq 0$ so $F'_{3}$ is not
slidable. The endpoints of $e^2_4$ in both vertices $v^2_1,v^2_2$ are now
located in the intervals $(e^2_3,e^2_1)$, so the endpoints of $e^1_4$ in
both vertices $v^1_1,v^1_2$ must also be located in the intervals
$(e^1_3,e^1_1)$; the endpoints of $e^1_4$ are indicated in
Fig.~\ref{n21c}(c) by open circles in $v^1_1,v^1_2$. The only possible
embedding of the edge $e^1_4$ in $T^1$ is shown in Fig.~\ref{n21c2}(a).

Now, by Lemma~\ref{mt} the bigons $F'_1,F'_2,F'_3$ form a complete disk
system of the handlebody $M_T$. Observe that the endpoints of $e^1_4$ and
$e^2_1$ can be connected via arcs $c_1,c_2$ along the strings
$I'_{1,2},I'_{2,1}$ that are disjoint from the corners of the bigons
$F'_1,F'_2,F'_3$. It follows that the circle $e^1_4\cup e^2_1\cup c_1\cup
c_2\subset\partial M_T$ is disjoint from $F'_1,F'_2,F'_3$ and hence
bounds a disk $F'_4$ in $M_T$ disjoint from $F'_1,F'_2,F'_3$. It is not
hard to see that the quotient $A=(F'_1\cup F'_2\cup F'_3\cup F'_4)/\psi$
is a surface in $(M,T_0)$ which intersects $T$ transversely and minimally
with $\Delta(\partial T,\partial A)=2$ and $G_{T,A}=T\cap A\subset T$ the
essential graph of Fig.~\ref{n21c2}(b). Since $T$ is positive, by the
parity rule $A$ must be a neutral annulus, hence each face of $G_{T,A}$
is necessarily a Scharlemann cycle of length 4, and since $M$ is
irreducible it follows that $A$ must separate $M$ into two solid tori
whose meridian disks, that is the faces of $G_{T,A}$, each intersect $A$
coherently in 4 arcs. Therefore $M$ is a Seifert fibered manifold over
the disk with two singular fibers of indices $4,4$, so $M(\partial T)$ is
a Seifert fibered torus bundle over the circle with horizontal bundle
fiber $\wh{T}$. By the classification of such torus bundles (cf
\cite[\S 2.2]{hatcher1}), it follows that $M(\partial
T)=(+0,0;1/2,-1/4,-1/4)$ and hence that $M=(+0,1;-1/4,-1/4)$.

\begin{figure}
\psfrag{g1}{$\gamma_1^1$}
\psfrag{g2}{$\gamma_2^2$}\psfrag{g}{$\gamma_1$}
\psfrag{e11}{$e^1_1$}\psfrag{e21}{$e^2_1$}
\psfrag{e12}{$e^1_2$}\psfrag{e22}{$e^2_2$}
\psfrag{e13}{$e^1_3$}\psfrag{e23}{$e^2_3$}
\psfrag{e14}{$e^1_4$}\psfrag{e24}{$e^2_4$}
\psfrag{e15}{$e^1_5$}\psfrag{e25}{$e^2_5$}
\psfrag{vj}{$w_i$}\psfrag{vj'}{$w_j$}
\psfrag{F'1}{$F'_1$}\psfrag{F'2}{$F'_2$}\psfrag{F'3}{$F'_3$}
\psfrag{F'4}{$F'_4$}\psfrag{F'5}{$F'_5$}\psfrag{F'6}{$F'_6$}
\psfrag{(a)}{$(a)$}\psfrag{(b)}{$(b)$}\psfrag{(c)}{$(c)$}
\psfrag{T1}{$T^1$}
\psfrag{T2}{$T^2$}\psfrag{T}{$T$}
\psfrag{e1}{$e_1$}\psfrag{e2}{$e_2$}
\psfrag{e3}{$e_3$}\psfrag{e4}{$e_4$}
\psfrag{e5}{$e_5$}\psfrag{e6}{$e_6$}
\psfrag{c1}{$c_1$}\psfrag{c2}{$c_2$}
\psfrag{c3}{$c_3$}\psfrag{c4}{$c_4$}
\Figw{5.5in}{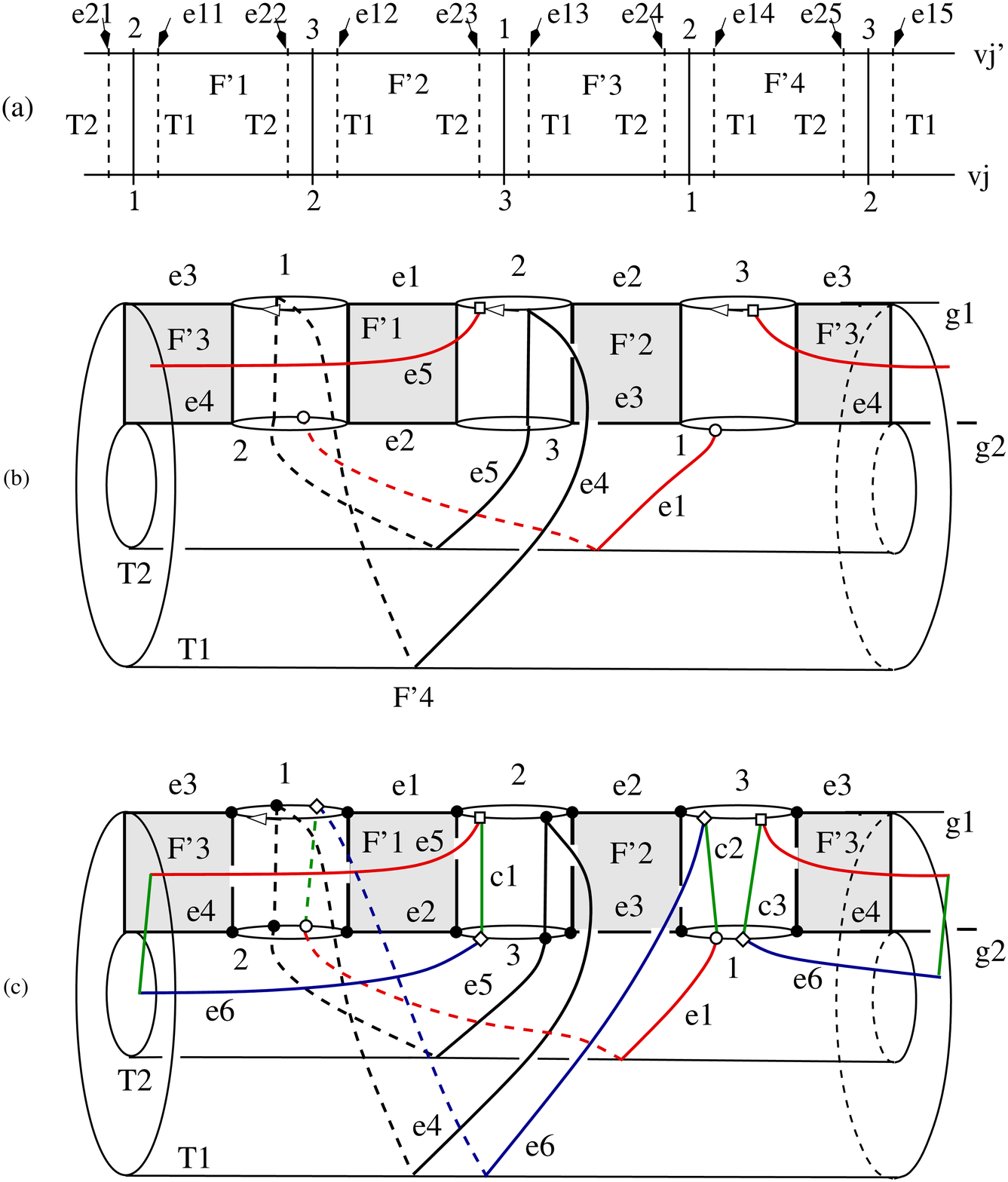}{}{n21c3}
\end{figure}

The case $t=3$ is handled in a similar way: Fig.~\ref{n21c3}(a) shows the
labeling of the edges $E=\{e_1,\dots,e_5\}$ and the bigon faces
$F'_1,\dots,F'_4$ they cobound in $G_R$, while Fig.~\ref{n21c3}(b) shows
the embedding of the bigons $F'_1,\dots,F'_4$ in $M_T$, up to
homeomorphism. Again there are only two possible embeddings for the edge
$e^2_1$ in $T^2$ depending on whether the endpoint of $e^2_1$ in $v^2_2$
lies in $(e^2_4,e^2_5)$ or $(e^2_5,e^2_2)$. In the first case the bigon
$F'_4$ is again slidable; we now sketch the details for the case
corresponding to the later option, where $F'_4$ is non-slidable.

Fig.~\ref{n21c3}(b) shows the embedding of $e^2_1\subset T^2$ which makes
the bigon $F'_4$ non-slidable and forces the embedding of $e^1_5\subset
T^1$. At this point all the edges in $E^2=\{e^2_1,\dots,e^2_5\}$ have
been embedded in $T^2$. Fig.~\ref{n21c3}(c) shows the embedding of a 6th
edge $e^2_6$ in $T^2$ with endpoints on the intervals
$(e^2_4,e^2_1)\subset v^2_1$ and $(e^2_4,e^2_1)\subset v^2_1$ and which
is disjoint from and not parallel in $T^2$ to any of the edges
$e^2_1,e^2_2,\dots,e^2_5$. These properties imply that
$e^1_6=\psi^{-1}(e^2_6)$ must be the edge in $T^1$ sketched in
Fig.~\ref{n21c3}(c).

\begin{figure}
\psfrag{T}{$T$}
\psfrag{e1}{$e_1$}\psfrag{e2}{$e_2$}
\psfrag{e3}{$e_3$}\psfrag{e4}{$e_4$}
\psfrag{e5}{$e_5$}\psfrag{e6}{$e_6$}
\Figw{4in}{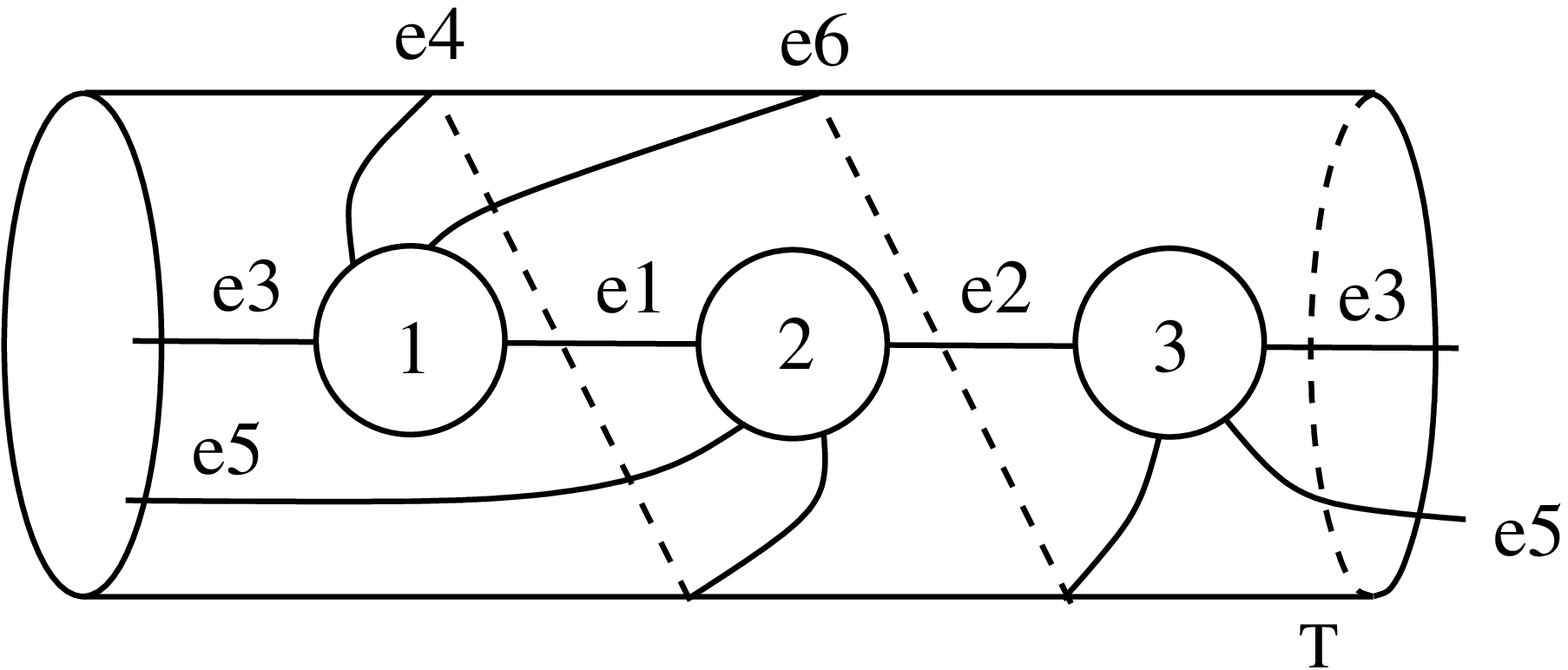}{}{n21c3b}
\end{figure}

As before, it is possible to connect the endpoints of $e^1_5$ and $e^2_6$
via corners $c_1\subset I'_{2,3}$ and $c_2\subset I'_{3,1}$, and the
endpoints of $e^1_6$ and $e^2_1$ via corners $c_3\subset I'_{3,1}$ and
$c_4\subset I'_{1,2}$, so that the circles $e^1_5\cup e^2_6\cup c_1\cup
c_2\subset\partial M_T$ and $e^1_6\cup e^2_1\cup c_3\cup
c_4\subset\partial M_T$ bound disks $F'_5\subset M_T$ and $F'_6\subset
M_T$, respectively, with the property that the enlarged collection of
disks $F'_1,\dots,F'_6$ is disjoint (see Fig.~\ref{n21c3}(c)).

It follows that $A=(F'_1\cup F'_2\cup \cdots\cup F'_6)/\psi$ is a neutral
annulus in $(M,T_0)$ which intersects $T$ transversely and minimally with
$G_{T,A}$ the essential graph of Fig.~\ref{n21c3b} and $\Delta(\partial
T,\partial A)=2$. Since $|\partial A|=2$, each face of $G_{T,A}$ is
necessarily a Scharlemann cycle and hence $A$ must separate $M$ into two
solid tori with meridian disks the faces of $G_{T,A}$. As $G_{T,A}$ has
two trigon faces on one side of $A$ and a 6-sided face in the other side
of $A$, we must have that $M$ is a Seifert manifold over the disk with
two singular fibers of indices $3,6$, which by the classification of such
torus bundles (cf \cite[\S 2.2]{hatcher1}) implies that $M(\partial
T)=(+0,0;1/2,-1/3,-1/6)$ and hence that $M=(+0,1;-1/3,-1/6)$. The lemma
follows.
\epf

\subsection{The generic case $t\geq 4$.}\label{t4}
Here we assume that $t\geq4$; we first show that for most values of
$\alpha$ the bigon $F'_{t+1}$ is slidable.

\begin{lem}\label{a1}
If $\alpha\not\equiv\pm 1\mod t$ then the bigon $F'_{t+1}$ is slidable.
\end{lem}

\bpf
If $t=4$ then $\alpha\equiv\pm 1\mod t$. For $t\geq 5$ the condition
$\alpha\not\equiv\pm 1\mod t$ is equivalent to saying that the pairs of
strings $\{I'_{1,2},I'_{1+\alpha,2+\alpha}\}$ and
$\{I'_{2,3},I'_{2+\alpha,3+\alpha}\}$ are disjoint. These four strings
are shown in Fig.~\ref{n05-3} along with the embeddings of the faces
$F'_1,\dots,F'_{t+1}$ from Fig.~\ref{n05-2}. As $e^2_{t+2}$ has one
endpoint on the interval $(e^2_{2-\alpha},e^2_2)\subset v^2_2$ and the
other on $(e^2_{2+\alpha},e^2_2)\subset v^2_{2+\alpha}$, the endpoints of
$e^1_{t+2}=\psi^{-1}(e^2_{t+2})$ must lie on
$(e^1_{2-\alpha},e^1_2)\subset v^1_2$ and $(e^1_{2+\alpha},e^1_2)\subset
v^1_{2+\alpha}$. Therefore the edge $e^1_{t+2}$ must be embedded in $T^1$
as shown in Fig.~\ref{n05-3}, which implies that the cycles $e_1\cup
e_{t+1}$ and $e_2\cup e_{t+2}$ have the same slope in $\wh{T}$ and hence
that $F'_{t+1}$ is slidable.
\epf

\begin{figure}
\psfrag{g1}{$\gamma_1^1$}
\psfrag{g2}{$\gamma_2^2$}
\psfrag{e22}{$e^2_2$}\psfrag{e11}{$e^1_1$}
\psfrag{e12}{$e^1_2$}\psfrag{e23}{$e^2_3$}
\psfrag{e22a}{$e^2_{2-\alpha}$}\psfrag{e22b}{$e^2_{2+\alpha}$}
\psfrag{e12a}{$e^1_{2-\alpha}$}\psfrag{e12b}{$e^1_{2+\alpha}$}
\psfrag{e2t2}{$e^2_{t+2}$}
\psfrag{e1t1}{$e^1_{t+1}$}
\psfrag{e1t2}{$e^1_{t+2}$}
\psfrag{T1}{$T^1$}
\psfrag{T2}{$T^2$}
\Figw{5in}{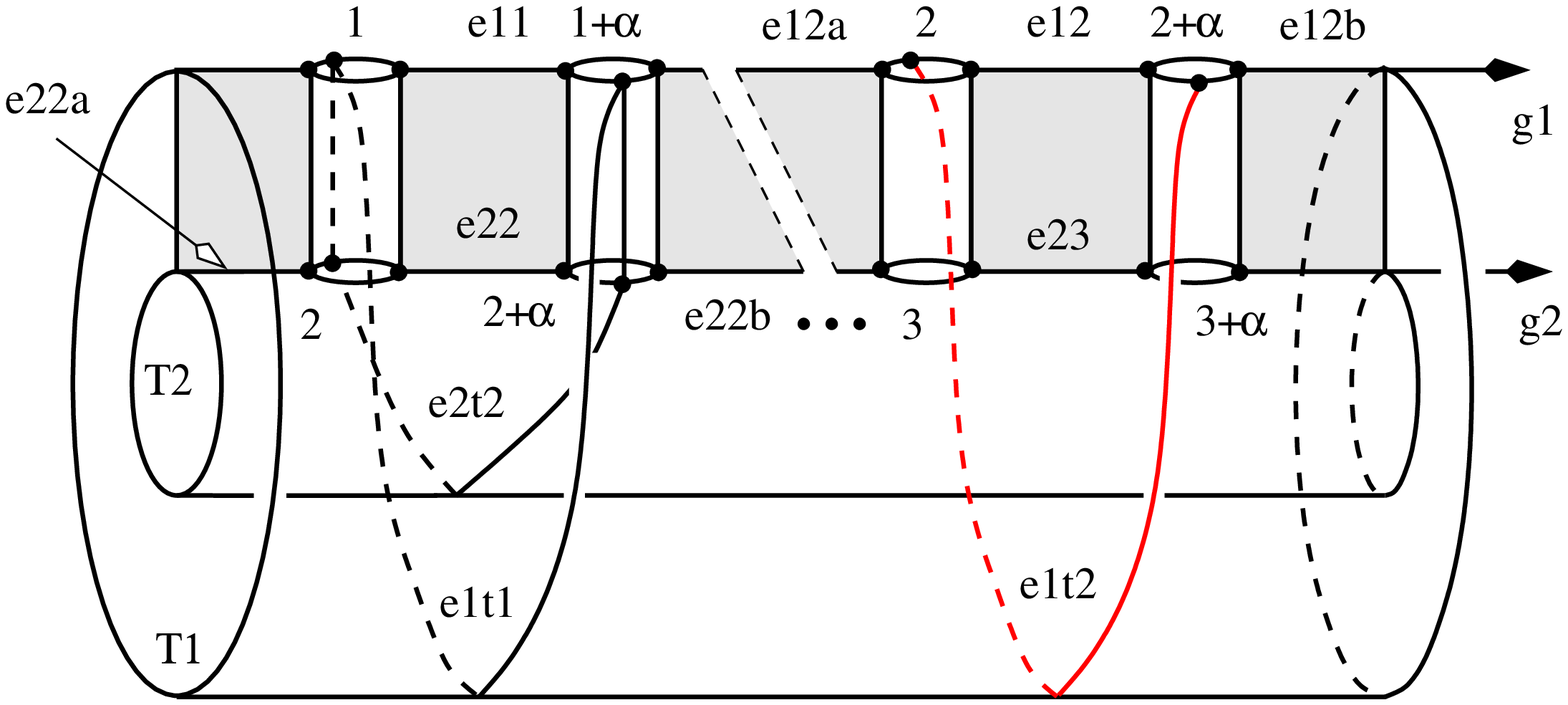}{}{n05-3}
\end{figure}

From here on we assume that $\alpha\equiv\pm 1\mod t$, and after
reversing the orientation of the edges of $E$, if necessary, we will take
$\alpha=1$.

\begin{lem}\label{alpha1}
If $\alpha=1$ and $t\geq 4$ then, up to homeomorphism, the bigons and
edges of $E$ are embedded in $M_T$ and $T$, respectively, as shown in
Fig.~\ref{n21-2a}(a),(b) if $F'_{t+1}$ is slidable, or
Fig.~\ref{n21-2b}(a),(b) if $F'_{t+1}$ is nonslidable. In the latter
case, if $\ove\subset\bg_R$ is the negative edge that contains the edges
of $E$ then $\ove=E$, so $|\ove|=t+2$.
\end{lem}

\bpf
By the argument of \cite[Lemma 3.6]{valdez11} (see in particular
\cite[Fig.\ 11]{valdez11} with $\alpha=1$) we may assume that, up to
homeomorphism, the bigons of $E$ are embedded in $M_T$ as shown in
Fig.~\ref{n21-2c}, where $t\geq 4$. Notice that the embedding of the
bigons of $E$ in $M_T$ only determines the embeddings of the edges
$e_1,e_2,\dots,e_{t+1}$ in $T^1$ and of $e_2,e_3,\dots,e_{t+2}$ in $T^2$.
To determine the possible embeddings of $e_1$ in $T^2$ and $e_{t+2}$ in
$T^1$ we proceed as follows.

\begin{figure}
\psfrag{g1}{$\gamma_1^1$}
\psfrag{g2}{$\gamma_2^2$}
\psfrag{F'1}{$F'_1$}\psfrag{F'2}{$F'_2$}\psfrag{F'3}{$F'_3$}
\psfrag{F'4}{$F'_4$}\psfrag{F't}{$F'_t$}
\psfrag{F't+1}{$F'_{t+1}$}\psfrag{F't-1}{$F'_{t-1}$}
\psfrag{T1}{$T^1$}
\psfrag{T2}{$T^2$}\psfrag{T}{$T$}
\psfrag{e1}{$e_1$}\psfrag{e2}{$e_2$}
\psfrag{e3}{$e_3$}\psfrag{e4}{$e_4$}
\psfrag{et}{$e_t$}\psfrag{et-1}{$e_{t-1}$}
\psfrag{et+1}{$e_{t+1}$}\psfrag{et+2}{$e_{t+2}$}
\Figw{6in}{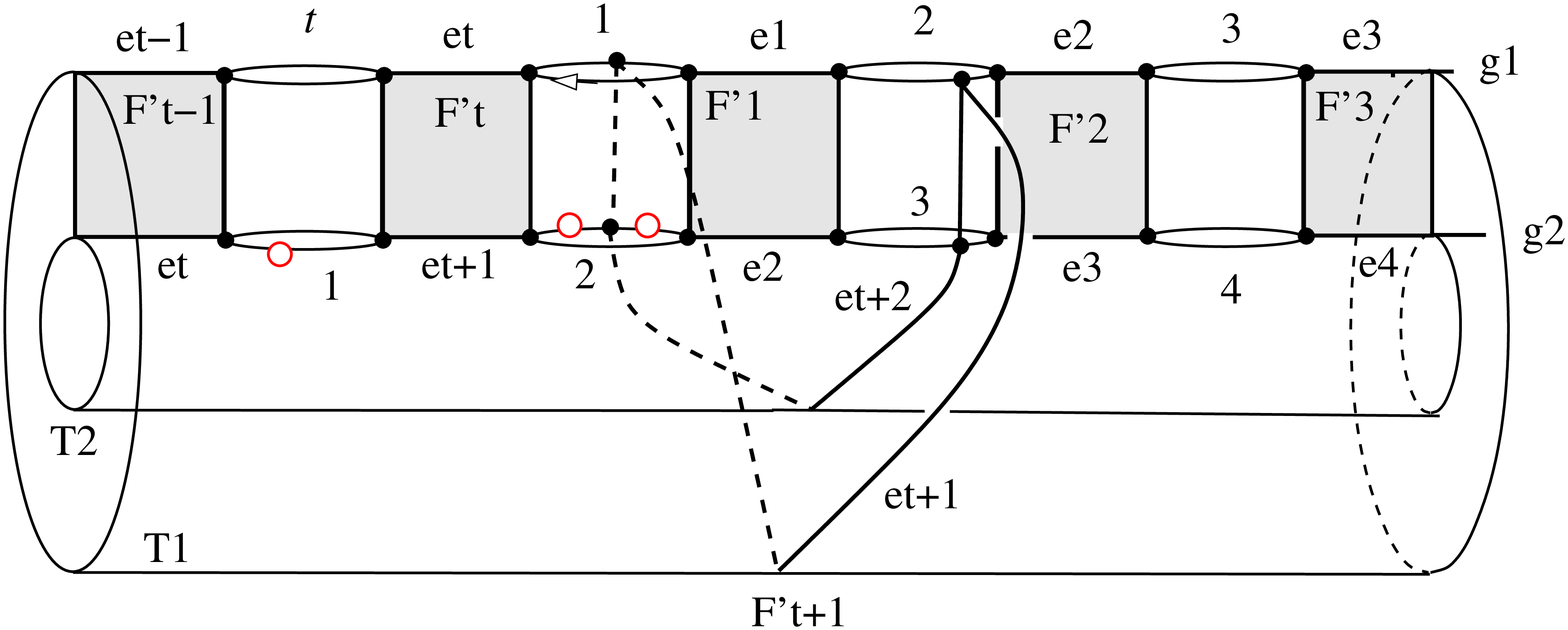}{}{n21-2c}
\end{figure}

From Fig.~\ref{n21-2c}, the endpoints of $e^1_1$ in $T^1$ are located in
the intervals $(e^1_{t+1},e^1_t)\subset v^1_1$ and
$(e^1_{t+1},e^1_2)\subset v^1_2$; since $\psi(T^1)=T^2$, it follows that
the endpoints of $e^2_1$ in $T^2$ must be located in the intervals
$(e^2_{t+1},e^2_t)\subset v^2_1$ and $(e^2_{t+1},e^2_2)\subset v^2_2$.

The interval $(e^2_{t+1},e^2_2)\subset v^2_2$ is split into the two
subintervals $(e^2_{t+1},e^2_{t+2})$ and $(e^2_{t+2},e^2_2)$ by the
endpoint of $e^2_{t+2}$ in $v^2_2$, which gives rise to two possible
locations of the endpoint of $e^2_1$ in $v^2_2$. The endpoints of $e^2_1$
are denoted by open circles in $v^2_1,v^2_2\subset T^2$ in
Fig.~\ref{n21-2c}. Since all faces of the graphs on $T^1,T^2$ of
Fig.~\ref{n21-2c} are disks, the embedding of $e^2_1$ in $T^2$ is
completely determined by the location of its endpoints; therefore there
are exactly two possible embeddings of $e^2_1$ in $T^2$.

\setcounter{case}{0}

\begin{case}
The endpoint of $e^2_1$ in $v^2_2$ lies in $(e^2_{t+1},e^2_{t+2})$.
\end{case}

This case is represented in Fig.~\ref{n21-2a}(a). Now the endpoints of
$e^2_{t+2}$ in $T^2$ lie in the intervals $(e^2_1,e^2_2)\subset v^2_2$
and $(e^2_3,e^2_2)\subset v^2_3$, and so the endpoints of
$e^1_{t+2}=\psi^{-1}(e^2_{t+2})$ in $T^1$ must lie in
$(e^1_1,e^1_2)\subset v^1_2$ and $(e^1_3,e^1_2)\subset v^1_3$. Therefore
$e^1_{t+2}$ must be the edge in $T^1$ shown in Fig.~\ref{n21-2a}(a). The
graph in $T$ produced by the edges of $E$ is shown in
Fig.~\ref{n21-2a}(b); as the cycles $e_1\cup e_{t+1}$ and $e_2\cup
e_{t+2}$ have the same slope in $\wh{T}$, the bigon $F'_{t+1}$ is
slidable.

\begin{figure}
\psfrag{g1}{$\gamma_1^1$}
\psfrag{g2}{$\gamma_2^2$}
\psfrag{F'1}{$F'_1$}\psfrag{F'2}{$F'_2$}\psfrag{F'3}{$F'_3$}
\psfrag{F'4}{$F'_4$}\psfrag{F't}{$F'_t$}
\psfrag{F't+1}{$F'_{t+1}$}\psfrag{F't-1}{$F'_{t-1}$}
\psfrag{(a)}{$(a)$}\psfrag{(b)}{$(b)$}
\psfrag{T1}{$T^1$}
\psfrag{T2}{$T^2$}\psfrag{T}{$T$}
\psfrag{e1}{$e_1$}\psfrag{e2}{$e_2$}
\psfrag{e3}{$e_3$}
\psfrag{e4}{$e_4$}
\psfrag{et}{$e_t$}
\psfrag{et-1}{$e_{t-1}$}
\psfrag{et+1}{$e_{t+1}$}\psfrag{et+2}{$e_{t+2}$}
\psfrag{etp1}{$e_{t+1}$}\psfrag{etp2}{$e_{t+2}$}
\Figw{6in}{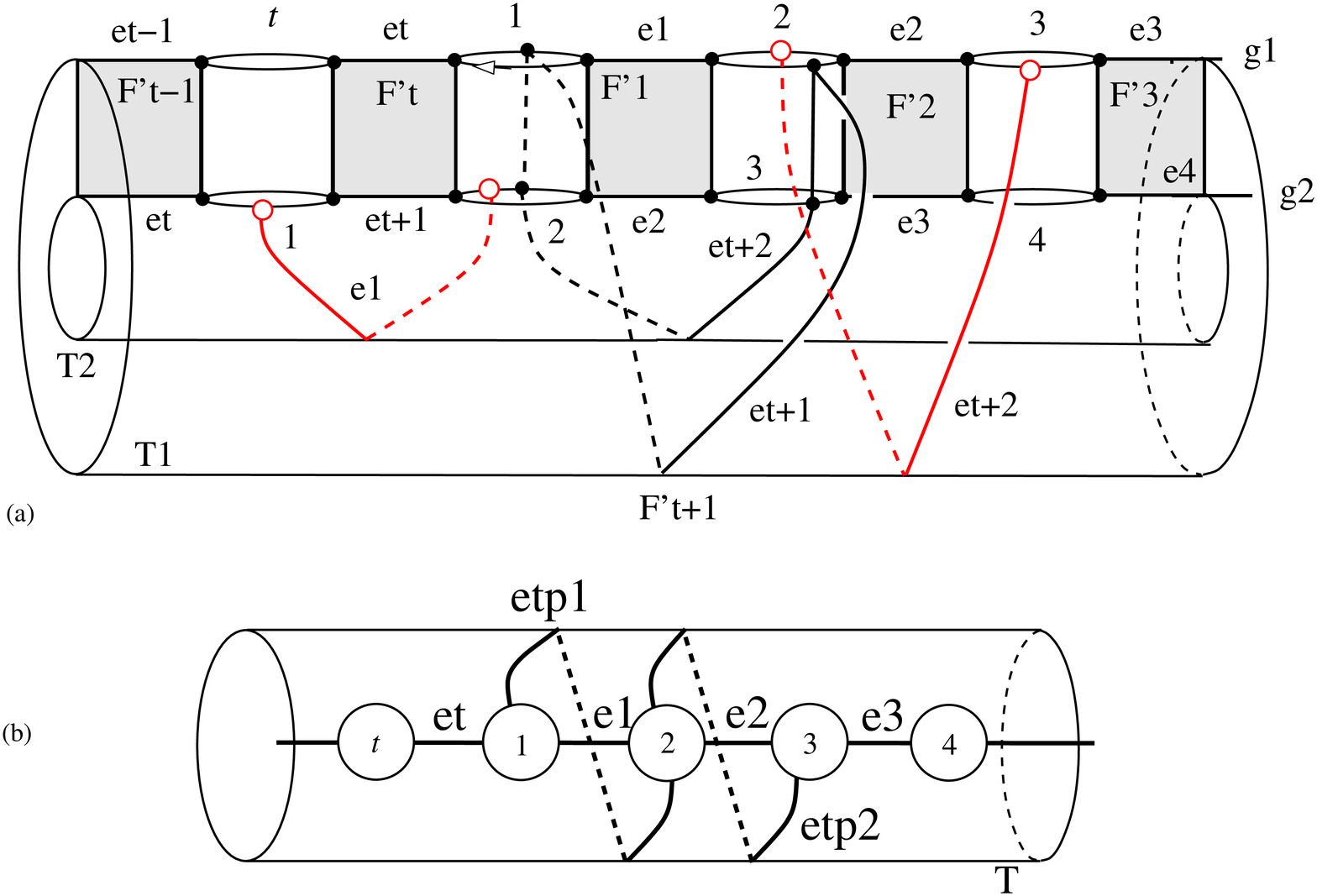}{}{n21-2a}
\end{figure}

\begin{case}
The endpoint of $e^2_1$ in $v^2_2$ lies in $(e^2_{t+2},e^2_{2})$.
\end{case}

The embedding of $e^2_1$ in $T^2$ is shown in Fig.~\ref{n21-2b}(a). The
endpoints of $e^2_{t+2}$ in $T^2$ lie in the intervals
$(e^2_{t+1},e^2_1)\subset v^2_2$ and $(e^2_3,e^2_2)\subset v^2_3$, so the
endpoints of $e^1_{t+2}$ in $T^1$ lie in $(e^1_{t+1},e^1_1)\subset v^1_2$
and $(e^1_3,e^1_2)\subset v^1_3$; thus $e^1_{t+2}$ must lie in $T^1$ as
shown in Fig.~\ref{n21-2b}(a). The graph in $T$ produced by the edges of
$E$ is shown in Fig.~\ref{n21-2b}(b); this time the cycles $e_1\cup
e_{t+1}$ and $e_2\cup e_{t+2}$ have slopes in $\wh{T}$ at distance 1 and
so the bigon $F'_{t+1}$ is nonslidable.

\begin{figure}
\psfrag{(a)}{$(a)$}\psfrag{(b)}{$(b)$}
\psfrag{g1}{$\gamma_1^1$}
\psfrag{g2}{$\gamma_2^2$}
\psfrag{F'1}{$F'_1$}\psfrag{F'2}{$F'_2$}\psfrag{F'3}{$F'_3$}
\psfrag{F'4}{$F'_4$}\psfrag{F't}{$F'_t$}
\psfrag{F't+1}{$F'_{t+1}$}\psfrag{F't-1}{$F'_{t-1}$}
\psfrag{T1}{$T^1$}
\psfrag{T2}{$T^2$}\psfrag{T}{$T$}
\psfrag{e1}{$e_1$}\psfrag{e2}{$e_2$}
\psfrag{e3}{$e_3$}
\psfrag{e4}{$e_4$}
\psfrag{et}{$e_t$}
\psfrag{et-1}{$e_{t-1}$}
\psfrag{et+1}{$e_{t+1}$}\psfrag{et+2}{$e_{t+2}$}
\psfrag{etp1}{$e_{t+1}$}\psfrag{etp2}{$e_{t+2}$}
\Figw{6in}{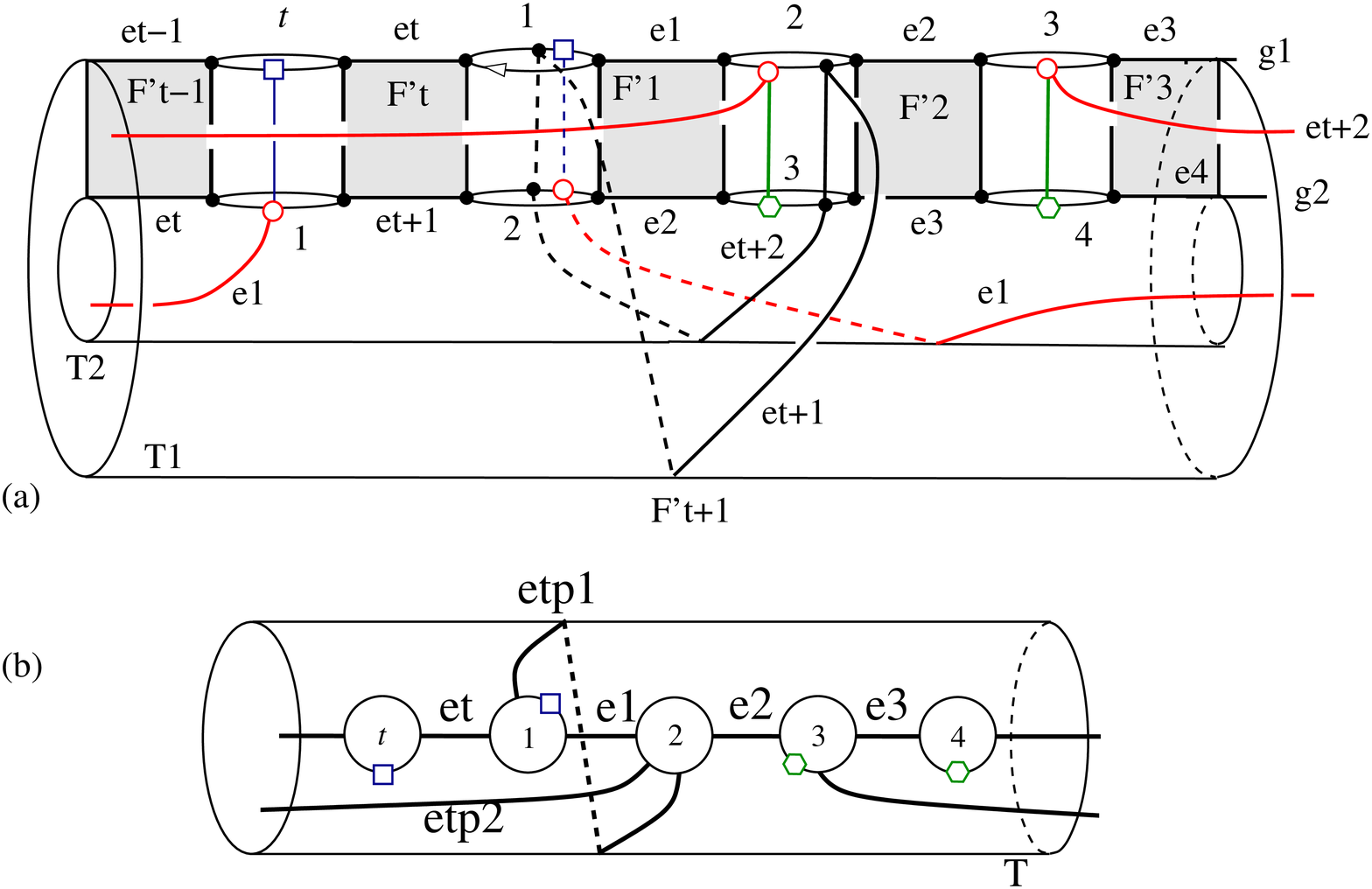}{}{n21-2b}
\end{figure}

Suppose there is a negative edge $e_0\subset\ove\setminus E$ in $G_R$
which cobounds a bigon $F_0$ with $e_1$, so that $F'_0\subset M_T$ is
cobounded by the edges $e^1_0\subset T^1$ and $e^2_1\subset T^2$. Since
$e^2_1$ has endpoints in $(e^2_{t+1},e^2_{t})\subset v^2_1$ and
$(e^2_{t+2},e^2_{2})\subset v^2_2$, the corners of $F'_0$ must lie in the
strings $I'_{t,1}$ and $I'_{1,2}$ and so the endpoints of $e^1_0$,
represented by squares in Fig.~\ref{n21-2b}, necessarily must lie in
$(e^1_{t},e^1_{t-1})\subset v^1_t$ and $(e^1_{t+1},e^1_{1})\subset
v^1_1$. This is impossible since such endpoints are separated by the edge
$e_{t+2}\subset T$ (see Fig.~\ref{n21-2b}(a),(b)), so no such edge $e_0$
exists in $\ove$. A similar argument shows that $\ove\setminus E$ does
not contain any edge adjacent to $e_{t+2}$. Therefore $E=\ove$ and so
$|\ove|=t+2$.
\epf

\begin{cor}\label{bd}
If $t\geq 4$ and $M$ is not homeomorphic to any of the manifolds $P\times
S^1/[m]$ then $|\ove|\leq t+2$ holds for any negative edge $\ove$ of
$\bg_R$.
\end{cor}

\bpf
Let $\ove$ be any negative edge of $\bg_R$ with $|\ove|\geq t+2$, and let
$E=\{e_1,\dots,e_{t+2}\}$ be any collection of $t+2$ consecutive edges in
$\ove$. Since $M$ is not any of the manifolds $P\times S^1/[m]$, by
Lemmas~\ref{lem01} and
\ref{a1} the permutation induced by $E$ must be of the form $x\to
x\pm 1$, and hence $|\ove|=t+2$ holds by Lemma~\ref{alpha1}.
\epf

We now prove Lemma~\ref{corr2}:

{\bf Proof of Lemma~\ref{corr2}: \ } If $t=3$ the bound $|\ove|\leq t+1$
holds for any negative edge $\ove$ of $\bgs$ by Lemma~\ref{t2t3}, so we
will assume that $t\geq 4$. Since none of the manifolds $P\times S^1/[m]$
is hyperbolic (cf
\cite[Proposition 3.4]{valdez11}), by Corollary~\ref{bd} the
bound $|\ove|\leq t+2$ holds for any negative edge $\ove$ of $\bgs$.

Suppose now there is an edge $\ove$ in $\bgs$ with $|\ove|=t+2$, so that
$T$ is a positive surface and hence, by the parity rule, in $G_S$ all
edges are negative, so any cycle in $G_S$ is even sided.

By \cite[Lemma 2.2(b)]{valdez11}, the reduced graph $\bgs$ has a vertex
$u$ of degree $n\leq 4$. Counting endpoints of edges of $G_S$ around $v$
yields the relations
$$
6t\leq\Delta\cdot t\leq n\cdot (t+2),
$$
which along with the restriction $n\leq 4$ imply that $n=t=4$,
$\Delta=6$, and that each of the 4 edges $\ove_1,\ove_2,\ove_3,\ove_4$
incident to $u$ in $\bgs$ has size $t+2=6$, as shown in Fig.~\ref{n30}.
We orient the edges of $\ove_1,\ove_2$ away from $u$ as indicated in
Fig.~\ref{n30}.

By Lemma~\ref{a1}, reversing the orientation of the vertices of $G_S$ and
relabeling the vertices of $G_T$, if necessary, we may assume that the
permutation induced by the oriented edge $\ove_1=\{e_1,\dots,e_6\}$
around $u$ is of the form $x\mapsto x+1$, while the permutation induced
by the oriented edge $\ove_2=\{a_1,\dots,a_6\}$ is of the form $x\mapsto
x\pm 1$.

\begin{figure}
\psfrag{e1}{$\ove_1$}\psfrag{e2}{$\ove_2$}
\psfrag{e3}{$\ove_3$}\psfrag{e4}{$\ove_4$}
\psfrag{a1}{$e_1$}\psfrag{a2}{$e_2$}\psfrag{a3}{$e_3$}
\psfrag{a4}{$e_4$}\psfrag{a5}{$e_5$}\psfrag{a6}{$e_6$}
\psfrag{b1}{$a_1$}\psfrag{b2}{$a_2$}\psfrag{b3}{$a_3$}
\psfrag{b4}{$a_4$}\psfrag{b5}{$a_5$}\psfrag{b6}{$a_6$}
\psfrag{v}{$u$}
\Figw{3in}{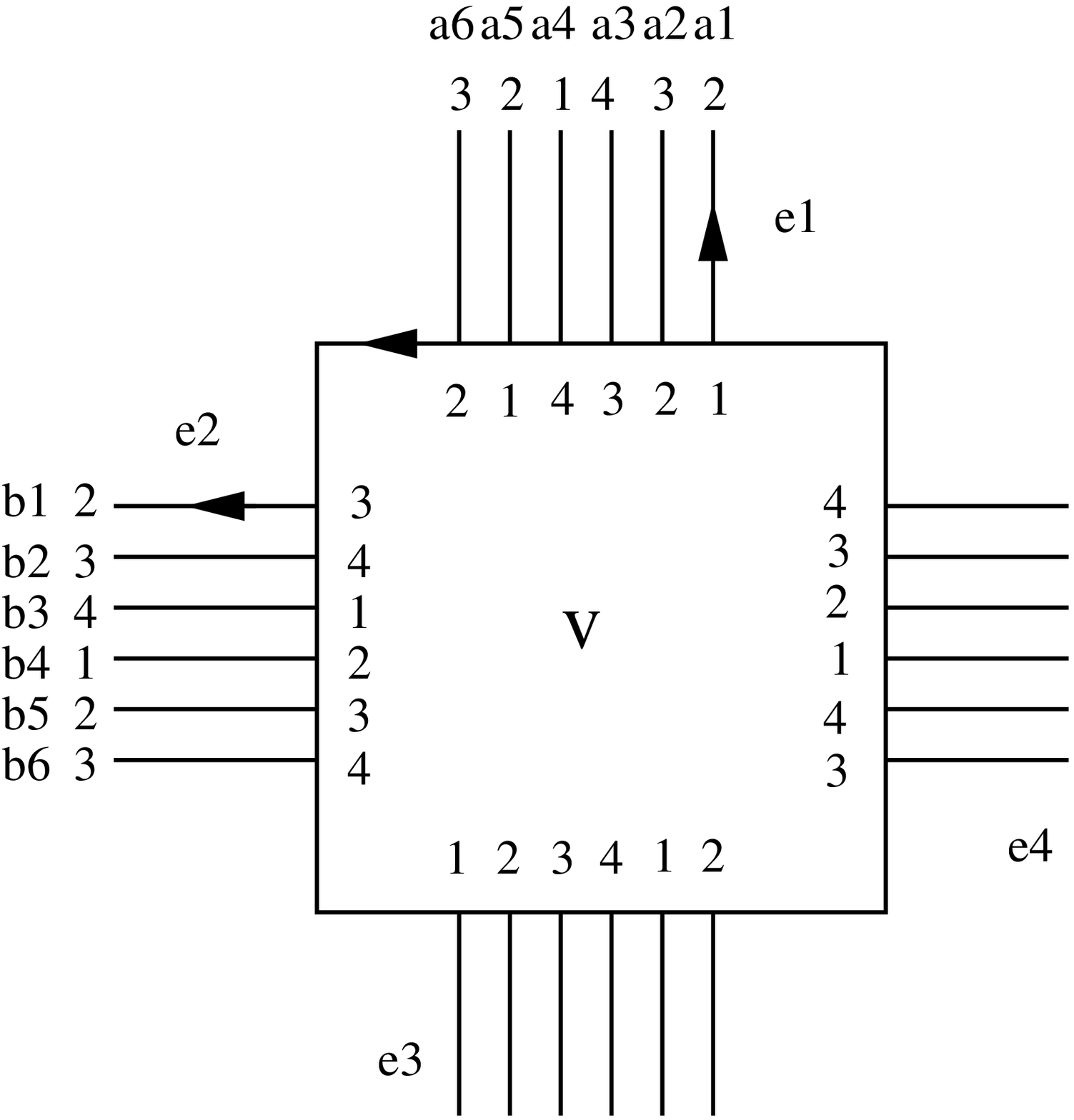}{}{n30}
\end{figure}

We consider in detail the case where $\ove_2$ induces the permutation
$x\mapsto x-1$; the case where this permutation is $x\mapsto x+1$ follows
by a similar argument. We may assume that the endpoints of the edges in
$\ove_1,\ove_2$ are labeled as in Fig.~\ref{n30}. We may also assume that
the edges of $\ove_1$ lie in $T$ as shown in Fig.~\ref{p01}; this figure
is obtained from the graph of Fig.~\ref{n21-2b}(b) with $t=4$, except
here $T$ is shown cut along the edges $e_1,e_5$ of $\ove_1$.

\begin{figure}
\psfrag{e1}{$\ove_1$}\psfrag{e2}{$\ove_2$}
\psfrag{e3}{$\ove_3$}\psfrag{e4}{$\ove_4$}
\psfrag{a1}{$e_1$}\psfrag{a2}{$e_2$}\psfrag{a3}{$e_3$}
\psfrag{a4}{$e_4$}\psfrag{a5}{$e_5$}\psfrag{a6}{$e_6$}
\psfrag{b1}{$a_1$}\psfrag{b2}{$a_2$}\psfrag{b3}{$a_3$}
\psfrag{b4}{$a_4$}\psfrag{b5}{$a_5$}\psfrag{b6}{$a_6$}
\psfrag{ai}{$a_i$}\psfrag{ai1}{$a_{i-1}$}
\psfrag{ai1?}{$a_{i-1}?$}
\psfrag{v}{$u$}\psfrag{T}{$T$}
\Figw{5in}{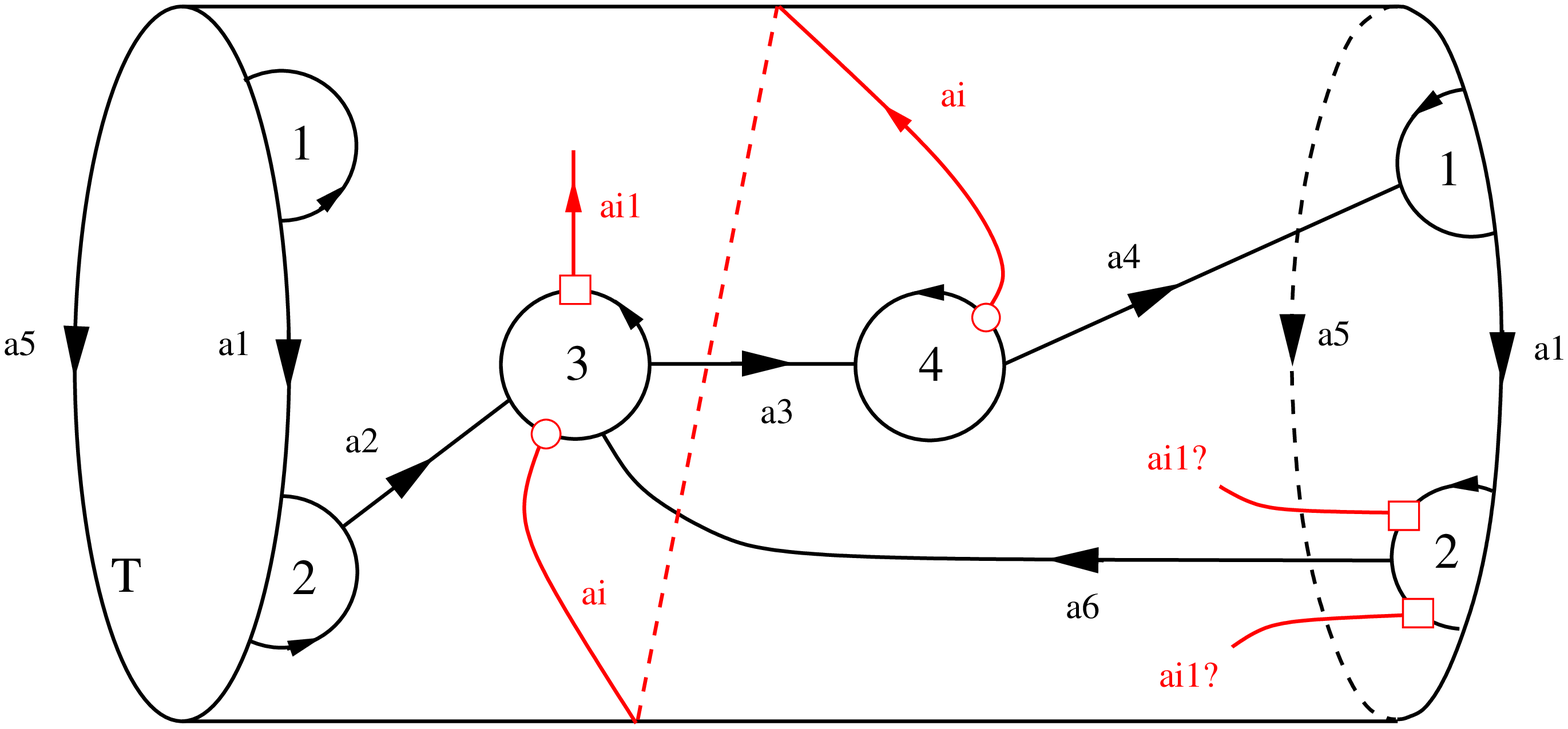}{}{p01}
\end{figure}

The edges $a_2,a_6$ of $\ove_2$ are parallel in $G_S$ and have endpoints
on the vertices $v_3,v_4$ of $G_T$, so by Lemma~\ref{basic}(b) for some
$i\in\{2,6\}$ the edge $a_i$ is not parallel in $T$ to $e_3$; it is not
hard to see that there is only one embedding of $a_i$ in $T$, as shown in
Fig.~\ref{p01}.

Let $F$ be the bigon of $\ove_2$ cobounded by $a_{i-1}$ and $a_i$; from
the point of view of $M_T$, the edge $a_{i-1}$ of $F$ lies in $T^1$ while
the edge $a_i$ lies in $T^2$. Since the edge $a^2_i\subset T^2$ has one
endpoint in the interval $(e^2_2,e^2_6)\subset v^2_3$ and its other
endpoint in the interval $(e^2_4,e^2_3)\subset v^2_4$, represented by
open circles in Fig.~\ref{p01}, by following the corners of $F$ up along
the strings $I'_{2,3}$ and $I'_{3,4}$ in Fig.~\ref{n21-2b}(a) (with
$t=4$) we can see that $a^1_{i-1}\subset T^1$ has one endpoint in the
interval $(e^1_1,e^1_5)\subset v^1_2$ and its other endpoint in the
interval $(e^1_3,e^1_2)$ of $v^1_3$. The possible locations of the
endpoints of $a_{i-1}$ around $v_2$ and $v_3$ are indicated by open
squares in Fig.~\ref{p01}; this situation is impossible since the edges
of $\ove_1\cup\{a_i\}\subset T$ separate the endpoints of $a_{i-1}$. This
contradiction shows that $|\ove|\leq t+1$ holds in $\bgs$.
\hfill\qed

\subsection{The manifolds $(M_t,T_0)$, $t\geq 4$}
For each $t\geq 4$ let $T$ be a $t$-punctured torus and $M_t$ the
manifold $M_T/\psi$, where $M_T=T\times I$ is the solid handlebody with
complete disk system $F'_1,F'_2,\dots,F'_{t+1}$ of Fig.~\ref{n21-2b}(a)
considered in Case 2 of Lemma~\ref{alpha1} and $\psi:T^1=T\times 0\to
T^2=T\times 1$ is the homeomorphism uniquely determined up to isotopy by
the conditions $\psi(v^1_i)=v^2_i$ for $1\leq i\leq t$ and
$\psi(e^1_j)=e^2_j$ for $1\leq j\leq t+1$. The basic properties of $M_t$
are summarized in the next lemma.

\begin{lem}\label{mt2}
The manifold $M_t$ is orientable and $\partial M_t$ is a torus $T_0$.
Moreover, the pair $(M_t,T_0)$ is irreducible and $T=T^1/\psi$ is a
properly embedded essential punctured torus in $(M_t,T_0)$.
\end{lem}

\bpf
If $T^1\subset\partial M_T$ and $T^2\subset\partial M_T$ are oriented by
the normal vectors $\vec{N}^1,\vec{N}^2$ in Fig.~\ref{n05-2} then the
conditions $\psi(v^1_i)=v^2_i$ for $1\leq i\leq t$ and
$\psi(e^1_j)=e^2_j$ for $1\leq j\leq t+1$ imply that the homeomorphism
$\psi:T^1\to T^2$ is orientation preserving, hence $M_t$ is orientable.
Clearly, $\partial M_t=(I'_{1,2}\sqcup I'_{2,3}\sqcup\cdots\sqcup
I'_{t,1})/\psi$ is a single torus $T_0$, and hence $T=T_1/\psi$ is a
properly embedded torus in $(M_t,T_0)$

If $T=T_1/\psi$ compresses in $M_t=M_T/\psi$ then, as $M_T$ is obtained
by cutting $M_t$ along $T$, it follows that $T^1$ or $T^2$ compresses in
$M_T$, which is not the case since $M_T=T\times I$. Therefore $T$ is
incompressible in $M_t$.

Since $M_T$ is a handlebody, hence irreducible, it now follows that $M_t$
is also irreducible. So, if $T_0$ compresses in $M_t$ then $M_t$ must be
a solid torus, contradicting the fact that $T$ is incompressible in
$M_t$. Therefore $T_0$ is incompressible and $T$ is essential in
$(M_t,T_0)$.
\epf

We now show that the bigons in the complete disk system of $M_T$ give
rise to a twice-punctured torus $S$ embedded in $(M_t,T_0)$.

\begin{lem}\label{d3}
For each $t\geq 4$ there is a separating, essential, twice-punctured
torus $S$ in $(M_t,T_0)$ such that $\Delta(\partial T,\partial S)=3$ and
$G_{S}=S\cap T\subset S$, $G_{T,S}=S\cap T\subset T$ are the graphs shown
in Figs.~\ref{n21-4} and \ref{n21-5}.
\end{lem}

\bpf
Let $t\geq 4$, and consider the faces $F'_1,\dots,F'_{t+1}$ and the edges
$e^i_j$ for $1\leq j\leq t+2$, $i=1,2$, embedded in the handlebody
$M_T=T\times I$ and $T^i$, respectively, as shown in
Fig.~\ref{n21-2b}(a), with $\psi(e^1_j)=e^2_j$ for all $j$. We add the
following objects to $M_T$ as indicated in Fig.~\ref{n21-3}:
\bit
\item
one bigon face parallel to $F'_i$ for $i=2,t$,

\item
two bigon faces parallel to $F'_i$ for $3\leq i\leq t-1$,

\item one more edge parallel to each of the edges
$e^2_{2},e^2_{3}$ and $e^1_{t},e^1_{t+1}$.
\eit

Therefore, in $T$, the edges $e_1,e_{t+2}$ get no parallel edges, each of
the edges $e_i$ for $i=2,t+1$ gets one parallel edge denoted by $e'_i$,
and each of the edges $e_i$ for $3\leq i\leq t$ gets two parallel edges
denoted $e'_i, e''_i$. The new collections of edges $e_i,e'_i,e''_i$
produce graphs $G^i\subset T^i$ for $i=1,2$ isomorphic to the graph in
Fig.~\ref{n21-4} such that, after a slight isotopy of the edges
$e'_i,e''_i$ if necessary, satisfy $\psi(G^1)=G^2$.

\begin{figure}
\psfrag{p1S}{$\partial_1 S$}\psfrag{p2S}{$\partial_2 S$}
\psfrag{g1}{$\gamma_1^1$}
\psfrag{g2}{$\gamma_2^2$}
\psfrag{F'1}{$F'_1$}\psfrag{F'2}{$F'_2$}\psfrag{F'3}{$F'_3$}
\psfrag{F'4}{$F'_4$}\psfrag{F't}{$F'_t$}
\psfrag{F't+1}{$F'_{t+1}$}\psfrag{F't-1}{$F'_{t-1}$}
\psfrag{T1}{$T^1$}
\psfrag{T2}{$T^2$}\psfrag{T}{$T$}
\psfrag{e1}{$e_1$}\psfrag{e2}{$e_2$}
\psfrag{e3}{$e_3$}\psfrag{e4}{$e_4$}
\psfrag{et}{$e_t$}\psfrag{et-1}{$e_{t-1}$}
\psfrag{et+1}{$e_{t+1}$}\psfrag{et+2}{$e_{t+2}$}
\psfrag{e'3}{$e'_3$}\psfrag{e''3}{$e''_3$}
\psfrag{e't+1}{$e'_{t+1}$}
\Figw{\the\textwidth}{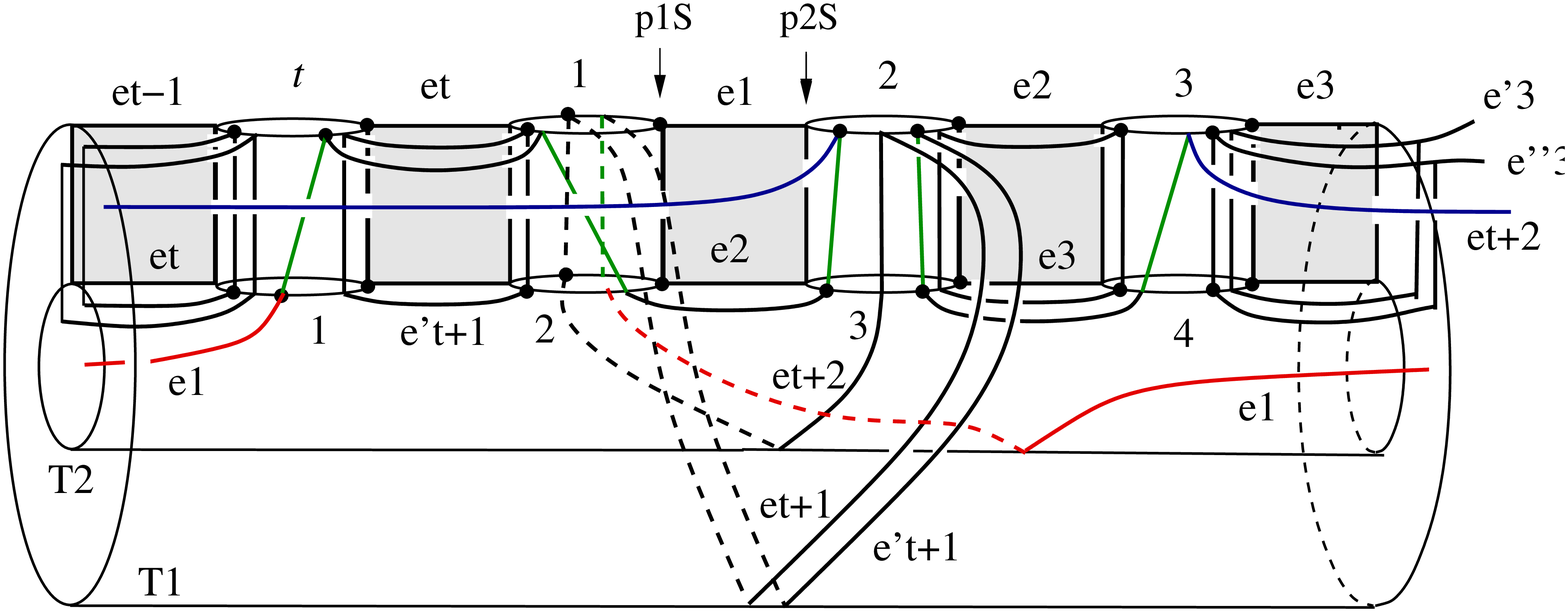}{}{n21-3}
\end{figure}

Now connect all the old and newly added edges of $T^1,T^2$ via mutually
disjoint corners as shown in Fig.~\ref{n21-3}. Observe that the 6-cycle
$\mc{C}$ in $\partial M_T$ containing the edges $e^2_1$, $(e'_{t+1})^1$,
$(e''_3)^2$, $e^1_{t+2}$, $(e'_{2})^2$, $(e''_t)^1$ is disjoint from the
complete disk system $F'_1,\dots,F'_{t+1}$ of $M_T$, hence $\mc{C}$
bounds a 6-sided disk `face' $F'_{\mc{C}}$ in $M_T$ disjoint from all the
other bigons in $M_T$.

\begin{figure}
\psfrag{T}{$T$}
\psfrag{e1}{$e_1$}\psfrag{e2}{$e_2$}
\psfrag{e3}{$e_3$}\psfrag{et}{$e_t$}
\psfrag{etp1}{$e_{t+1}$}\psfrag{etp2}{$e_{t+2}$}
\Figw{4.5in}{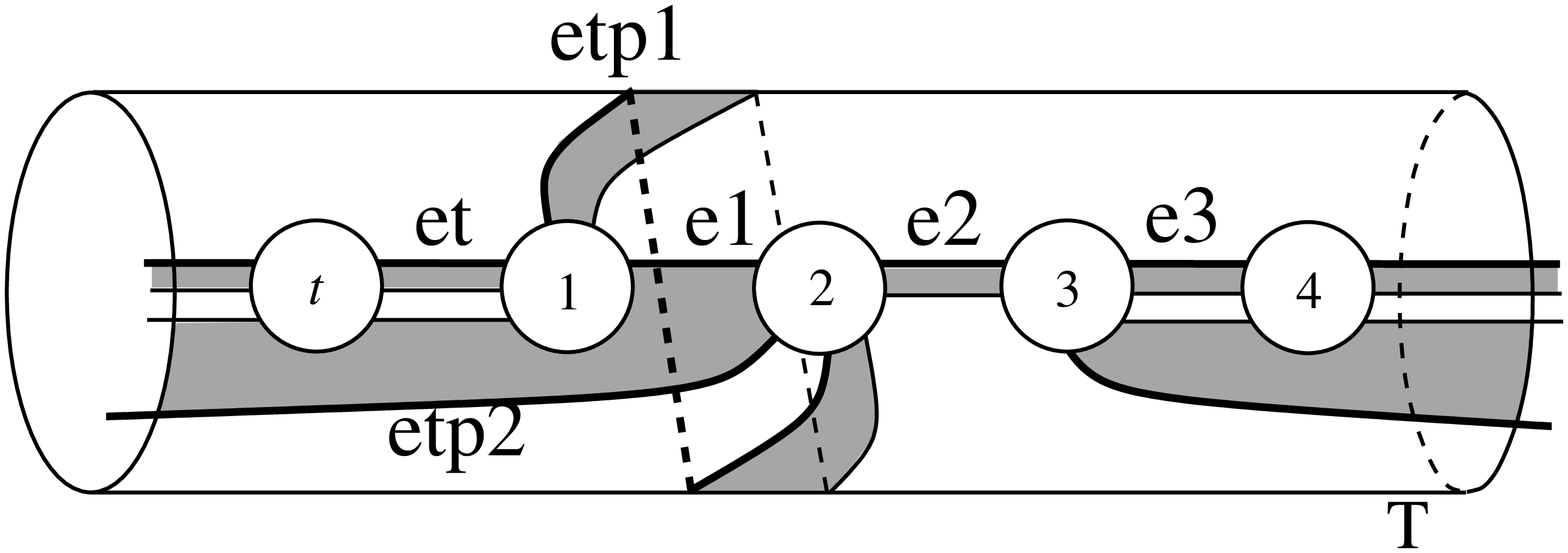}{}{n21-4}
\end{figure}

In this way we obtain a collection $\mc{F}$ of $3t-2$ disjoint disk faces
embedded in $M_T$: $2(3)+3(t-3)=3t-3$ bigons and one 6-sided disk face.
The condition $\psi(G^1)=G^2$ guarantees that $S=\mc{F}/\psi\subset
(M,T_0)$ is a properly embedded surface which intersects $T$ transversely
in the graph $G_{T,S}=T\cap S\subset T$ of Fig.~\ref{n21-4}. Moreover,
the collection of corners of faces in $\mc{F}$ whose endpoints are capped
with a small closed disk in Fig.~\ref{n21-3} form one boundary component
$\partial_1 S$ of $S$, while the remaining corners form a second boundary
component $\partial_2 S$; thus $|\partial S|=2$.

Therefore all faces of the graph $G_S=S\cap T\subset S$ are disks, and
$G_S$ has 2 vertices, $3t$ edges, and $|\mc{F}|=3t-2$ faces, so $\wh{S}$
is a surface with Euler number 0; since each vertex of $G_{T,S}$ has
degree 6, it follows that $\Delta(\partial S,\partial T)=3$.

Now, the faces of $G_{T,S}$ can be colored black or white as shown in
Fig.~\ref{n21-4}, and this coloring induces a corresponding black and
white coloring of the components of $M_T\setminus\cup\mc{F}$ such that
each face in $\mc{F}$ is adjacent to opposite colored components, which
implies that $S$ separates $M_t$ and hence that $S$ is a 2-punctured
separating torus. We denote by $S^B,S^W$ the closures of the components
of $M_t\setminus S$.

The graph $G_S$ can be determined by following the boundary circle
$\partial_1 S$ in Fig.~\ref{n21-3} in the direction of increasing labels;
starting at the endpoint of $e^1_1$ in $v_1\subset\partial T$ we find
that the collections of edges $\{e_1,\dots,e_{t+2}\}$,
$\{e'_2,\dots,e'_{t+1}\}$, and $\{e''_3,\dots,e''_{t}\}$ form distinct
parallelism classes in $T$ and are read in this order as we traverse the
circle $\partial_1 S$ in Fig.~\ref{n21-3}, and so $G_{S}$ must be the
graph shown in Fig.~\ref{n21-5}.

For each $*\in\{B,W\}$ the manifold $S^*$ is irreducible with $\partial
S^*$ a genus two surface. Notice that $G_{T,S}$ contains Scharlemann
cycle disk faces in $S^*$ of distinct lengths; any such Scharlemann cycle
disk face is nonseparating in $S^*$, whence $S^*$ is a genus two
handlebody with complete disk system any pair of Scharlemann cycle disk
faces of $G_{T,S}$ in $S^*$ of different sizes.

For $S^B$ we can take as complete disk system the bigon $x$ of $G_{T,S}$
containing the edge $e_2$ and the $t$-sided face $y$ containing
$e_{t+2}$, so that $\pi_1(S^B)=\set{x,y
\ | \ -}$ and $\partial_1S\subset\partial S^B$ is represented
by the word in $\pi_1(S^B)$ obtained by reading the consecutive
intersections of $\partial_1 S$ with the disks $x$ and $y$. Following
$\partial_1 S$ around in Fig.~\ref{n21-3}, we can see that
$\partial_1S=(yx)^2y^{t-2}$, which is not a primitive word in
$\pi_1(S^B)=\set{x,y \ | \ -}$; therefore $S$ is incompressible in $S^B$
by \cite[Lemma 5.2]{valdez7}. Similarly, in $S^W$ we can take as complete
disk system the white bigon $X$ in $G_{T,S}$ with corners along $v_3$ and
$v_4$ and the white $t+4$-sided face $Y$ containing $e_{t+2}$; we then
compute that $\partial_1S\subset\partial S^W$ is represented by the word
$Y^{t+3}XYX$ in $\pi_1(S^W)=\set{X,Y
\ | \ -}$, which is not
primitive, so $S$ is incompressible in $S^W$ too. Therefore $S$ is
incompressible, hence essential, in $M_t$.
\epf

\begin{figure}
\psfrag{S}{$S$}
\psfrag{ei}{$e_i$}\psfrag{e'i}{$e'_i$}
\psfrag{e''i}{$e''_i$}
\psfrag{F'1}{$F'_1$}\psfrag{F'C}{$F'_{\mc{C}}$}
\psfrag{p1S}{$\partial_1 S$}\psfrag{p2S}{$\partial_2 S$}
\Figw{3.5in}{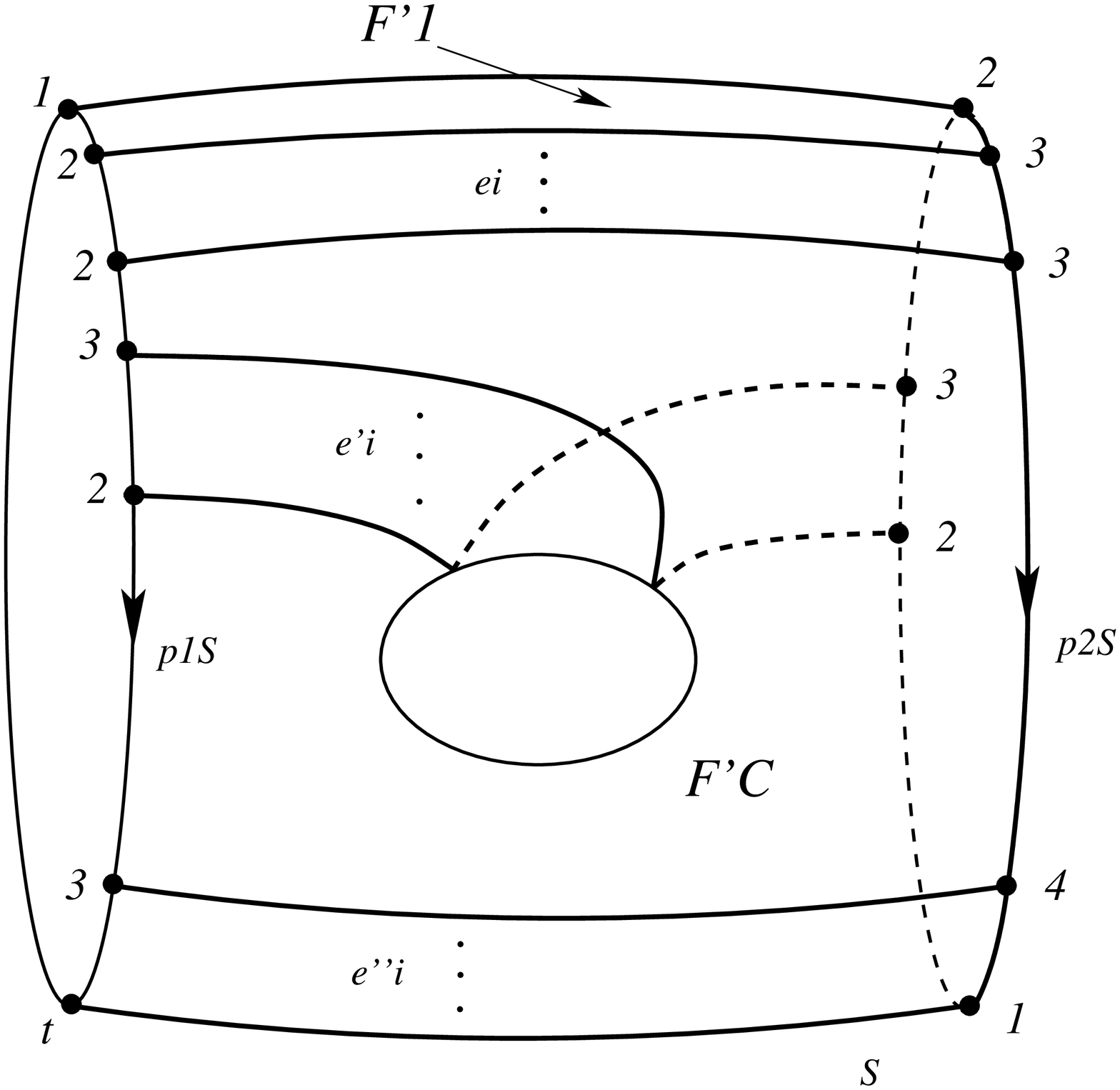}{}{n21-5}
\end{figure}

We now complete the proofs of the remaining main results of this paper.

{\bf Proof of Proposition~\ref{main}: \ } Parts (1), (2) and (3) follow
immediately from Lemmas~\ref{mt2} and
\ref{d3}. By Lemma~\ref{mt}, the manifold $M_t(\partial T)$ is a
torus bundle over a circle with fiber $\wh{T}$. Finally, by
\cite[Proposition 3.4]{valdez11}(c), if $(P\times S^1/[m],T_0)$
contains two $\mc{K}$-incompressible tori $T,T'$ then
$\Delta(T,T')\in\{0,1,2,4\}$; since $M_t$ contains the essential torus
$S$ with $\Delta(\partial T,\partial S)=3$, it follows that $M_t$ is not
homeomorphic to any of the manifolds $P\times S^1/[m]$, so part (4)
holds.
\hfill\qed

{\bf Proof of Proposition~\ref{prop2}: \ } That $T$ is positive follows
from Lemma~\ref{basic}(c). Suppose $(M,T_0)$ does not satisfy parts (1)
and (3) of the proposition. By Lemmas~\ref{lem01}(1) and \ref{t2t3} we
then have that $t\geq 4$, and so by Lemmas~\ref{lem01}(2), \ref{a1},
\ref{alpha1}, and the definition of $M_t$, we have that
$(M,T_0)$ is homeomorphic to $(M_t,T_0)$. Since, by
Proposition~\ref{main}, $M_t$ is not a manifold of the form $P\times
I/[m]$, the bound $|\ove|\leq t+2$ holds for all negative edges $\ove$ of
$\bg_R$ by Corollary~\ref{bd}. Therefore part (2) holds.
\hfill\qed

{\bf Proof of Lemma~\ref{tp2a}: \ } Assume parts (1), (2) and (4) of the
lemma do not hold; then $t\geq 4$ and $(M,T_0)\approx (M_t,T_0)$ by
Proposition~\ref{prop2}(2), and there is a negative edge in $\bg_Q$ of
length $|\ove|\geq t+2$. By Corollary~\ref{bd} we then have that
$|\ove|=t+2$, so the lemma holds.
\hfill\qed


\end{document}